%% file: hho-vi.tex
\documentclass[a4paper,11pt]{article}
\usepackage[english]{babel}
\usepackage[utf8]{inputenc}
\usepackage[T1]{fontenc}
\usepackage{amsmath,amsfonts,amssymb,mathtools}
\DeclareSymbolFontAlphabet{\mathbb}{AMSb}
\DeclareMathAlphabet{\mathbbmsl}{U}{bbm}{m}{sl}
\usepackage{nicefrac}
\usepackage{enumitem}
\usepackage{graphicx,placeins}
\usepackage[dvipsnames, svgnames]{xcolor}
\usepackage{hyperref}
\hypersetup{
 colorlinks = true,
 linkcolor = OliveGreen,
 citecolor = violet,
 filecolor = blue,
 urlcolor = magenta
}
\usepackage[hmargin={28mm,28mm},vmargin={30mm,35mm}]{geometry}
\usepackage{amsthm}
\usepackage{authblk}
\usepackage{newtxtext,newtxmath}
\usepackage[makeroom]{cancel}
\allowdisplaybreaks
\usepackage{pgfplots,pgfplotstable,booktabs}
\usepgfplotslibrary{groupplots}
\usetikzlibrary{calc}
\pgfplotsset{compat=1.18}
\newtheorem{theorem}{Theorem}
\newtheorem{lemma}[theorem]{Lemma}

\newtheorem{corollary}[theorem]{Corollary}
\theoremstyle{remark}
\newtheorem{remark}[theorem]{Remark}
\theoremstyle{definition}
\newtheorem{definition}[theorem]{Definition}

\newcommand{\email}[1]{\href{mailto:#1}{\textcolor{blue}{#1}}}

\newcommand{\defi}{\mathrel{\mathop:}=}

\newcommand{\st}{~:~}

\renewcommand{\vec}[1]{\boldsymbol{#1}}
\newcommand{\normal}{\vec{n}}

\newcommand{\Norm}[2]{\left\|#2\right\|_{#1}}

\newcommand{\tnorm}[2]{|\kern-0.25ex|\kern-0.25ex|#2|\kern-0.25ex|\kern-0.25ex|_{#1}}

\newcommand{\Natu}{\mathbb{N}}
\newcommand{\Real}{\mathbb{R}}
\newcommand{\Comp}{\mathbb{C}}
\newcommand{\ii}{\mathrm{i}}

\DeclareMathOperator{\Div}{div}

\DeclareMathOperator{\Grad}{grad}

\newcommand{\LL}[1][\Omega]{L^2(#1)}
\newcommand{\bLL}[1][\Omega]{\vec{L}^2(#1)}

\newcommand{\Hzdiv}[1][\Omega]{\vec{H}_{0}(\Div ; #1)}

\newcommand{\Poly}{\mathcal{P}}
\newcommand{\bPoly}{\vec{\Poly}}

\renewcommand{\t}{\mathcal{T}}
\newcommand{\f}{\mathcal{F}}

\newcommand{\proj}[1]{\pi^{\ell}_{#1}}
\newcommand{\vproj}[1]{\vec\pi^{\ell}_{#1}}

\newcommand{\djdiv}{\underline{\mathbbmsl{J}}_h^{\ell}}
\newcommand{\djdivT}{\underline{\mathbbmsl{J}}_T^{\ell}}
\newcommand{\djdivz}{\underline{\mathbbmsl{J}}_{h,0}^{\ell}}

\newcommand{\jh}{\underline{\mathbbmsl{j}}_h}
\newcommand{\jTF}{\underline{\mathbbmsl{j}}_T}
\newcommand{\jT}{\vec{j}_T}
\newcommand{\jFn}{j_{F,n}}
\newcommand{\jTh}{\vec{j}_h}

\newcommand{\wh}{\underline{\mathbbmsl{w}}_h}
\newcommand{\wTF}{\underline{\mathbbmsl{w}}_T}
\newcommand{\wT}{\vec{w}_T}
\newcommand{\wFn}{w_{F,n}}
\newcommand{\wTh}{\vec{w}_h}

\newcommand{\DpT}{D_T^{\ell}}
\newcommand{\Dph}{D_h^{\ell}}
\newcommand{\IhT}{\underline{\mathbbmsl{I}}^{\ell}_T}
\newcommand{\Ihh}{\underline{\mathbbmsl{I}}^{\ell}_h}

\newcommand{\LLstar}{L^2_{\star}(\Omega)}
\newcommand{\Ph}{\Poly^{\ell}_{\star}(\t_h)}

\usepackage[normalem]{ulem}
\newcounter{corr}
\definecolor{violet}{rgb}{0.580,0.,0.827}
\newcommand{\corr}[3]{\typeout{Warning : a correction remains in page \thepage}
	\stepcounter{corr}
	{\color{blue}\ifmmode\text{\,\sout{\ensuremath{#1}}\,}\else\sout{#1}\fi}
	{\color{red}#2}
	{\color{violet} #3}
}

\begin{document}

\title{A hybrid high-order discretization of~a~volume~integral~formulation for eddy currents on~unstructured~polyhedral~meshes}
\author[1]{Raman Kumar}
\author[1]{Aurelio Spadotto}
\author[1,2]{Silvano Pitassi}
\affil[1]{IMAG, Univ Montpellier, CNRS, Montpellier, France,
\email{raman.kumar@umontpellier.fr},
\email{aurelio-edoardo.spadotto@umontpellier.fr}}
\affil[2]{Graduate School of Mathematical Sciences, The University of Tokyo, Tokyo, Japan,
\email{silvano-pitassi@g.ecc.u-tokyo.ac.jp}}
\date{}
\maketitle

\begin{abstract}
Volume integral formulations are attractive for eddy-current simulation because they require a mesh only of the conducting region. However, their dense interaction matrices make it important to control the number of unknowns needed for a prescribed accuracy. We propose an arbitrary-order Hybrid High-Order discretization of a current-based volume integral formulation on general unstructured polyhedral meshes. The method combines high-order polynomial approximation with local mesh refinement around defects and material interfaces. For bounded positive spatially varying resistivity and magnetic permeability uniform throughout space, we establish well-posedness and a priori error estimates using a constructive right inverse of the discrete divergence. Numerical experiments demonstrate convergence for smooth solutions and compare mesh refinement with increasing the polynomial degree. Selected comparisons show that higher degree on coarser meshes can reduce the number of current unknowns coupled by the Green operator and the associated storage for a prescribed accuracy. Finally, a three-dimensional eddy-current testing benchmark assesses the predicted impedance against experimental measurements, using local refinement around a narrow slot and reusing the system matrix and its factorization throughout the probe scan.
\par\medskip
\noindent\textbf{Keywords:}
Volume integral equations; Eddy currents; Hybrid high-order methods;
Polyhedral meshes; Mixed methods; Eddy current testing.
\par
\end{abstract}

\section{Introduction}
\label{sec:introduction}
\begin{figure}[t]
\centering
\includegraphics[width=0.7\textwidth]{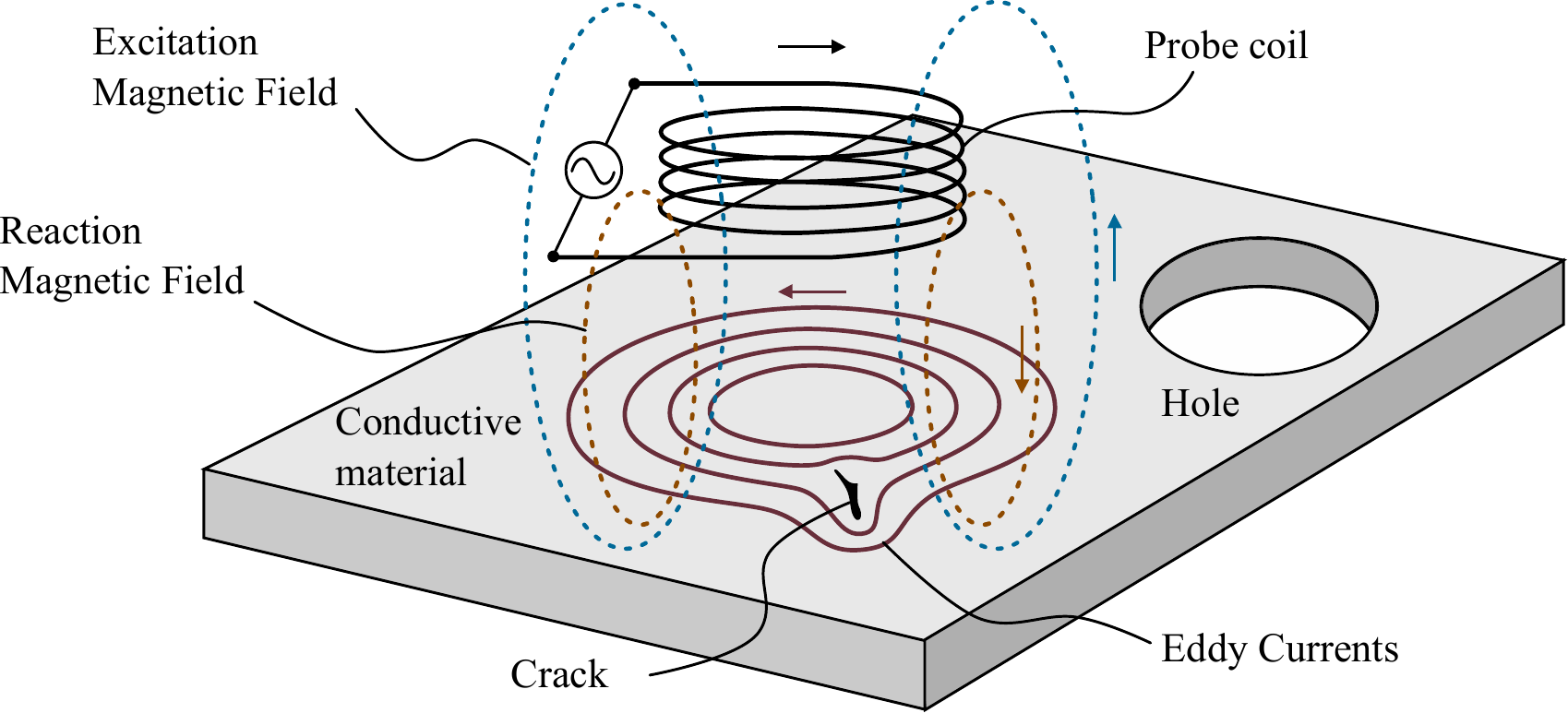}
\caption{Schematic configuration of an eddy-current testing problem. A probe coil generates an excitation magnetic field that induces eddy currents inside the conducting material. In the presence of a localized defect, such as a crack or a hole, the current distribution is perturbed, which in turn modifies the reaction magnetic field measured by the probe. This type of configuration motivates the use of volume integral formulations.}
\label{fig:ect-configuration}
\end{figure}

Volume integral formulations for eddy currents express the electromagnetic reaction in terms of currents supported in the conductor, so the surrounding air need not be meshed~\cite{Passarotto.Pitassi.Specogna:22}.
This feature is particularly useful in \emph{eddy-current testing (ECT)}, where a probe scans a conducting component containing a localized defect.
The defect changes the current paths and the coil impedance, while the conductor geometry remains fixed; see Figure~\ref{fig:ect-configuration}.
For an air-cored probe in a magnetically homogeneous environment, the same conductor mesh and system matrix can be used throughout a scan at fixed frequency and material parameters. Only the excitation and impedance evaluation change with the probe position.
ECT is the main application motivating the methodology presented here.

\noindent How the defect enters the model depends on its physical nature. A local change in conductivity is represented through the resistivity coefficient,
whereas material loss, such as a hole or an insulating slot, is represented by the boundary of the conducting domain.
Interpreting an impedance scan requires understanding how the signal depends on defect shape and local resistivity.

\noindent The computational cost of this conductor-only description is a dense Green interaction operator.
Its dimension depends on the current approximation space, making the choice of spatial resolution central even when only the conductor is discretized.
High-order approximation is of particular interest as a means of attaining a prescribed accuracy with fewer globally interacting current unknowns.
General polyhedral meshes provide additional geometric flexibility.
Their faces can represent material interfaces and defect boundaries, small cells can resolve their neighbourhoods, and larger cells can cover regions that admit coarser resolution.
The accuracy and computational cost obtained by combining mesh resolution and polynomial degree are a central concern of this work.

\noindent We develop an arbitrary-order \emph{Hybrid High-Order (HHO)} discretization of a mixed, current-density-based volume integral formulation for three-dimensional eddy currents on general polyhedral meshes.
To the best of our knowledge, this is the first HHO discretization of three-dimensional eddy-current problems accompanied by a stability and convergence analysis and an experimental benchmark validation.
The continuous formulation, introduced in Section~\ref{sec:continuous.model}, retains the current density as its primary unknown and provides a direct starting point for this study.
In the discrete scheme of Section~\ref{sec:hho.scheme}, the current is represented by cell polynomials and single-valued normal face unknowns, while a broken scalar multiplier enforces a reconstructed divergence constraint.
The Green operator acts directly on the cell polynomials, combining an explicit representation for nonlocal integration with a hybrid conservation structure.

\noindent Three-dimensional volume integral eddy-current methods have been developed using edge- and facet-based representations, including current approximations adapted to the solenoidality constraint~\cite{Nguyen.Meunier.ea:14}.
Low-order compatible schemes on general polyhedral meshes, including discrete geometric approaches closely related to mimetic finite differences, were introduced in~\cite{Bettini.Passarotto.Specogna:17}.
Iterative realizations avoiding explicit construction of the fully populated inductance matrix were considered in~\cite{Bettini.Passarotto.Specogna:18}.
The low-order polyhedral framework was systematized in~\cite{Passarotto.Pitassi.Specogna:22}, with formal links to mimetic reconstruction principles in~\cite{Pitassi:21}.
Optimized cycle bases~\cite{Voltolina.Torchio.ea:21} and low-rank compression~\cite{Vacalebre.Pitassi.Specogna:23} address the cost of nonlocal coupling.
The present method extends the approximation choices available within a polyhedral current-based integral formulation to arbitrary polynomial degree.

\noindent A different approach uses a discontinuous Galerkin discretization on polyhedral meshes, analysed within the framework of partially coercive Friedrichs systems~\cite{Torosyan.Bannwarth.ea:25}.
Its recent tensor-train realization~\cite{Torosyan.Bannwarth.ea:26} includes three-dimensional benchmarks and compresses the integral operator after tensorization.
The hybrid construction considered here uses shared normal face unknowns and a commuting divergence reconstruction to impose discrete current conservation.
Higher-order spline methods provide further context, including cohomology representations for eddy-current formulations~\cite{Kapidani.Merkel.ea:24} and an isogeometric augmented electric field integral equation~\cite{Nolte.Torchio.ea:25}.
Our approximation spaces are defined directly on the cells and faces of an unstructured polyhedral mesh and do not require a volumetric patch parametrization.

\noindent The local mixed structure is closely related to arbitrary-order mixed methods for diffusion~\cite{Di-Pietro.Ern:17} and HHO discretizations of Darcy flow~\cite{Chave.Di-Pietro.Formaggia:18,Di-Pietro.Droniou:20}.
Commuting divergence operators and hybrid flux unknowns already play a central role in that literature.
Here, boundedness of the nonlocal term in $L^2$ and coercivity of the resistive term preserve the underlying mixed stability mechanism.
The integral formulation nevertheless changes the algebraic structure and integration.
The Green operator couples cell currents across the conductor and involves a weakly singular kernel.
Section~\ref{sec:implementation} describes this structure, the resulting Green block dimension, and the quadrature procedure, including the consequences for the usual HHO static condensation.

\noindent We provide a complete stability and error analysis, whose main results are stated in Section~\ref{sec:main.results} and proved in Section~\ref{sec:stability.error}.
The inf-sup proof constructs a right inverse of the discrete divergence by combining a relative cellular cochain lifting with a local polynomial completion.
This gives an alternative to the Fortin argument used, for example, in~\cite{Chave.Di-Pietro.Formaggia:18}, and a concrete instance of the discrete Poincar\'e framework for HHO differential forms developed in~\cite{Pitassi:26}.
Related hybrid stability constructions for div--curl systems and discrete Weber inequalities appear in~\cite{Dalphin.Ducreux.Lemaire.Pitassi:25,Lemaire.Pitassi:25}.
The lifting controls the divergence constraint, while boundedness of the Green operator enters the continuity and consistency estimates, yielding well-posedness and a priori error estimates.
The estimates allow cellwise regularity across fitted material interfaces and transfer directly to bounded linear observations of the current.
A componentwise normalization of the scalar potential extends the construction to finitely many separated conductors, with their magnetic interaction retained by the global Green operator.

\noindent The numerical study in Section~\ref{sec:numerical.results} verifies the approximation properties and assesses the method as a computational tool.
A manufactured solution on Cartesian, tetrahedral, and Voronoi meshes tests convergence.
Comparisons between finer degree-one and coarser degree-two discretizations quantify accuracy against the dimension and storage of the Green block, together with assembly and solution costs.
An experimental ECT benchmark with a fastener hole and a thin radial slot then tests the impedance response on locally resolved polyhedral meshes and separates the fixed matrix setup from repeated solves along the probe scan.
The computations use an uncompressed Green operator, providing a baseline for combining mesh and degree choices with fast integral-operator techniques.
Section~\ref{sec:conclusions} discusses these computational developments and extensions to other integral formulations.

\section{Continuous model}
\label{sec:continuous.model}
Let $\Omega$ be a bounded, connected Lipschitz domain in $\Real^3$ representing the conductor.
The domain need not be simply connected and may contain holes and finite-width slots, as in the ECT benchmark below. We consider linear time-harmonic eddy currents with angular frequency $\omega>0$ and use the phasor convention $\Re(\widehat{\vec v}\,\mathrm e^{\ii\omega t})$.

\noindent The magnetic permeability $\mu>0$ is uniform throughout space. The resistivity $\rho\in L^\infty(\Omega)$ is real-valued and satisfies $0<\rho_\ast\leq\rho\leq\rho^\ast$ almost everywhere in $\Omega$. It may therefore vary within the conductor and jump across material interfaces.
A jump in $\mu$, even across an interface aligned with mesh faces, would require an additional magnetic response in the integral representation and is outside the present model.
An inclusion with large but finite resistivity remains part of $\Omega$, whereas an insulating cavity or slot is excluded from $\Omega$ and carries a zero normal-current condition on its walls. The aluminium benchmark uses constant resistivity and $\mu=\mu_0$.

\noindent All function spaces are understood over $\Comp$, and $L^2$ inner products are linear in their first argument and conjugate-linear in their second argument.

\noindent We consider a current-density-based volume integral formulation~\cite{Bossavit:81,Passarotto.Pitassi.Specogna:22}. Given a magnetic vector potential $\vec{a}_0\in\bLL$, find the current density $\vec{j}$ and the electric scalar potential $u$ such that
\begin{equation}
\label{eq:strong.model}
\left\{
\begin{aligned}
\rho\vec{j}+\ii\omega\mu G(\vec{j})+\Grad u &= -\ii\omega\vec{a}_0 &&\text{in }\Omega,\\
\Div\vec{j} &= 0 &&\text{in }\Omega,\\
\vec{j}\cdot\normal &= 0 &&\text{on }\partial\Omega.
\end{aligned}
\right.
\end{equation}
In the first equation, $\rho\vec j$ is the resistive electric field, $\mu G(\vec j)$ is the induced magnetic vector potential, and $-\ii\omega\vec a_0$ is the electric excitation associated with the prescribed source potential.
The scalar potential $u$ enforces current conservation together with the insulating boundary condition.
The Green potential operator $G$ is defined by
\[
G(\vec{j})(\vec{x}) \coloneqq\int_{\Omega} g(\vec{x},\vec{y})\,\vec{j}(\vec{y})\,\mathrm{d}\vec{y}, \qquad
\text{with} \qquad g(\vec{x},\vec{y})\defi \frac{1}{4\pi\lvert \vec{x}-\vec{y}\rvert}.
\]

\noindent Schur's test gives boundedness of the Newtonian operator $G$ on $\bLL$, which we use below to estimate the nonlocal term. Since the kernel is real and symmetric, $G$ is self-adjoint and $\ii\omega\mu(G\vec v,\vec v)_\Omega$ has zero real part.
The quadratic form also has a magnetic energy interpretation~\cite{Passarotto.Pitassi.Specogna:22}. Let $\vec A$ be the Newtonian potential of the zero extension of $\vec v\in\bLL$. Then $-\Delta\vec A$ equals that extension, and integration by parts gives
\[
  (G\vec v,\vec v)_\Omega = \sum_{i=1}^3
       \Norm{L^2(\Real^3)}{\nabla A_i}^2 \geq 0.
\]

\noindent For the scalar unknown, we impose the zero-average condition by introducing the space
\[
\LLstar \defi \left\{ q\in\LL\mid \int_{\Omega} q =0\right\}.
\]
The mixed weak formulation is to find $(\vec{j},u)\in \Hzdiv\times\LLstar$ such that
\begin{subequations}
\label{eq:weak.model}
\begin{align}
(\rho\vec{j},\vec{w})_{\Omega} +\ii\omega(\mu G(\vec{j}),\vec{w})_{\Omega}-(u,\Div\vec{w})_{\Omega}
&=(-\ii\omega\vec{a}_0,\vec{w})_{\Omega} &&\forall \vec{w}\in\Hzdiv,
\label{eq:weak.model.a}\\
(\Div\vec{j},q)_{\Omega}&=0 &&\forall q\in\LLstar.
\label{eq:weak.model.b}
\end{align}
\end{subequations}
\noindent This formulation keeps the divergence-free condition in mixed form and does not require a current vector potential or an explicit construction of topological cycles.
Connectedness is used to fix the scalar potential by a single zero-average condition.
For a conductor with several connected components, the corresponding condition is imposed on each component, while the Green operator still couples their currents.

\section{HHO scheme}
\label{sec:hho.scheme}
We first define the hybrid current and scalar spaces, then introduce the reconstructed divergence, stabilization, and discrete mixed problem.
The reduction operators relate these constructions to the continuous fields and establish their commuting and polynomial-consistency properties.
The section closes with the main stability and error estimates.

\subsection{Broken polynomial spaces and hybrid unknowns}
\label{sec:discrete.spaces}

Let $(\t_h)_h$ be a family of conforming polyhedral meshes of $\Omega$, with planar faces collected in $\f_h$.
For each cell $T\in\t_h$, let $\f_T$ be its set of faces, and denote the cell and face diameters by $h_T$ and $h_F$, respectively; set $h=\max_{T\in\t_h}h_T$.
We assume mesh regularity in the usual HHO sense~\cite{Di-Pietro.Droniou:20}, including uniformly shape-regular conforming simplicial submeshes $(\mathcal S_h)_h$ with uniformly bounded local complexity.
In particular, for each tetrahedron $S\in\mathcal S_h$ contained in $T$ and each two-dimensional simplex $\widehat F$ contained in $F$, uniformly in $h$,
\[
h_S\simeq h_T,
\qquad
h_{\widehat F}\simeq h_F.
\]

\noindent We fix a unit normal $\normal_F$ on each face $F\in\f_h$, and, for every $T\in\t_h$ and $F\in\f_T$, we set
\[
\varepsilon_{T,F}\defi \normal_{T,F}\cdot\normal_F\in\{-1,1\},
\]
where $\normal_{T,F}$ denotes the unit normal to $F$ pointing outward from $T$.

\noindent Given an integer $\ell\in\Natu$, we introduce the (broken) polynomial spaces
\[
\Poly^{\ell}(\t_h)
\defi
\left\{
q_h\in\LL
\mid
q_{h\mid T}\in\Poly^{\ell}(T)
\quad\forall T\in\t_h
\right\},
\]
and
\[
\bPoly^{\ell}(\t_h)
\defi
\left\{
\vec v_h\in\bLL
\mid
\vec v_{h\mid T}\in\bPoly^{\ell}(T)
\quad\forall T\in\t_h
\right\}.
\]
We denote by $\pi_T^{\ell}$ and $\pi_F^{\ell}$ the $L^2$-orthogonal projectors onto $\Poly^{\ell}(T)$ and $\Poly^{\ell}(F)$, respectively, and by $\vec\pi_T^{\ell}$ the $\bLL[T]$-orthogonal projector onto $\bPoly^{\ell}(T)$.

\noindent We define the $\vec H(\Div)$-like hybrid space
\begin{equation}
\label{eq:djdiv}
\djdiv
\coloneqq
\left\{
\jh
\defi
\Bigl(
(\jT)_{T\in\t_h},
(\jFn)_{F\in\f_h}
\Bigr)
\st
\begin{alignedat}{2}
\jT&\in\bPoly^{\ell}(T)
&&\qquad\forall T\in\t_h,
\\
\jFn&\in\Poly^{\ell}(F)
&&\qquad\forall F\in\f_h
\end{alignedat}
\right\}.
\end{equation}
For a generic element $\jh\in\djdiv$ and each cell $T\in\t_h$, we denote its local restriction by
\[
\jTF
\defi
\bigl(
\jT,
(\jFn)_{F\in\f_T}
\bigr)
\in\djdivT.
\]
We  associate with $\jh$ the broken cellwise polynomial field $\jTh\in\bPoly^{\ell}(\t_h)$ defined by
\[
(\jTh)_{\mid T}\defi\jT
\qquad
\forall T\in\t_h.
\]
On boundary faces, we define the discrete normal trace operator
\begin{equation}
\label{eq:discrete.trace}
\gamma_{n,h}:\djdiv\to L^2(\partial\Omega) \;\; \text{ and } \;\;
\gamma_{n,h}(\jh)_{\mid F}\defi\jFn
\qquad
\forall F\in\f_h\cap\partial\Omega.
\end{equation}
The discrete counterpart of $\Hzdiv$ is the subspace of $\djdiv$ defined by
\begin{equation}
\label{eq:djdivz}
\djdivz
\defi
\left\{
\jh\in\djdiv
\mid
\gamma_{n,h}(\jh)=0
\right\}.
\end{equation}
For the scalar multiplier, we consider the broken polynomial space with zero average
\begin{equation}
\label{eq:Ph}
\Ph
\defi
\left\{
q_h\in\Poly^{\ell}(\t_h)
\mid
(q_h,1)_\Omega=0
\right\}.
\end{equation}

\subsection{Divergence reconstruction and stabilization}
\label{sec:reconstruction}

The discrete divergence is reconstructed cellwise through the standard HHO integration-by-parts formula~\cite{Di-Pietro.Droniou:20}, involving the cell polynomial and the single-valued normal face unknowns.

\begin{definition}[Local divergence reconstruction]
For each $T\in\t_h$, the local divergence reconstruction $\DpT:\djdivT\to\Poly^{\ell}(T)$ is defined by requiring, for every $\jTF\in\djdivT$,
\begin{equation}
\label{eq:def.DT}
\bigl(\DpT(\jTF),q\bigr)_T
=
-(\jT,\Grad q)_T
+
\sum_{F\in\f_T}
\varepsilon_{T,F}
(\jFn,q_{\mid F})_F
\qquad
\forall q\in\Poly^{\ell}(T).
\end{equation}
\end{definition}
\noindent The associated global reconstruction $\Dph:\djdiv\to\Poly^{\ell}(\t_h)$ is defined cellwise by
\[
(\Dph(\jh))_{\mid T}
\defi
\DpT(\jTF)
\qquad
\forall T\in\t_h.
\]
Notice that
\begin{equation}
\label{eq:Dph.zero.mean}
\Dph(\djdivz)\subset\Ph.
\end{equation}
Indeed, taking $q=1$ in~\eqref{eq:def.DT} and summing over the mesh, the internal face contributions cancel by single-valuedness of the face unknowns and the orientation convention, while the boundary contributions vanish for functions in $\djdivz$.

\noindent For all $\wTF,\jTF\in\djdivT$, we define the local stabilization by
\begin{equation}
\label{eq:stab.div}
s_{\Div,T}(\jTF,\wTF)
\coloneqq
\sum_{F\in\f_T}
h_F
\Bigl(
\vec j_{T\mid F}\cdot\normal_F-\jFn,
\vec w_{T\mid F}\cdot\normal_F-\wFn
\Bigr)_F.
\end{equation}
The corresponding global stabilization is
\begin{equation}
\label{eq:global.stab.div}
s_h(\jh,\wh)
\defi
\sum_{T\in\t_h}
s_{\Div,T}(\jTF,\wTF).
\end{equation}
The weight $h_F$ gives the stabilization the same $L^2$ scaling as the cell-current mass term.

\subsection{Discrete bilinear forms}
\label{sec:bilinear.forms}
For all $\jh,\wh\in\djdivz$ and $q_h\in\Ph$, we discretize the resistive term using the cell unknowns as
\begin{equation}
\label{eq:def.mh}
m_h(\jh,\wh)
\defi
(\rho\jTh,\wTh)_\Omega.
\end{equation}
Since the Green operator enters without derivatives, it is applied directly to the broken cell polynomial field, yielding
\begin{equation}
\label{eq:def.gh}
g_h(\jh,\wh)
\defi
(\mu G(\jTh),\wTh)_\Omega.
\end{equation}
Although nonlocal and giving rise to a dense matrix, this term involves cell unknowns only.

\noindent The divergence coupling is defined by
\begin{equation}
\label{eq:def.bh}
b_h(\wh,q_h)
\defi
-
\sum_{T\in\t_h}
\bigl(
\DpT(\wTF),
q_{h\mid T}
\bigr)_T
=
-(\Dph\wh,q_h)_\Omega.
\end{equation}
Under the convention adopted above, $b_h$ is linear in its first argument and conjugate-linear in its second argument.

\noindent We then introduce the main sesquilinear form
\begin{equation}
\label{eq:def.ah}
a_h(\jh,\wh)
\defi
m_h(\jh,\wh)
+
\ii\omega g_h(\jh,\wh)
+
\alpha s_h(\jh,\wh),
\end{equation}
where $\alpha>0$ is a fixed stabilization parameter.
The analysis holds for any fixed $\alpha>0$.

\noindent Finally, the discrete right-hand side is given by
\begin{equation}
\label{eq:def.lh}
\ell_h(\wh)
\defi
(-\ii\omega\vec a_0,\wTh)_\Omega.
\end{equation}
The HHO discretization is formulated as the following mixed problem.
\begin{definition}[HHO scheme]
Find $(\jh,u_h)\in\djdivz\times\Ph$ such that
\begin{subequations}
\label{eq:discrete.problem}
\begin{align}
a_h(\jh,\wh)
+
\overline{b_h(\wh,u_h)}
&=
\ell_h(\wh)
&&
\forall\wh\in\djdivz,
\label{eq:discrete.problem.a}
\\
b_h(\jh,q_h)
&=
0
&&
\forall q_h\in\Ph.
\label{eq:discrete.problem.b}
\end{align}
\end{subequations}
\end{definition}
\noindent Upon inserting the definitions of the discrete forms, the first equation in \eqref{eq:discrete.problem} can equivalently be written as
\begin{multline}
\label{eq:expanded.scheme}
(\rho\jTh,\wTh)_\Omega
+
\ii\omega(\mu G(\jTh),\wTh)_\Omega
+
\alpha
\sum_{T\in\t_h}
s_{\Div,T}(\jTF,\wTF)
\\
-
\sum_{T\in\t_h}
\bigl(
u_{h\mid T},\DpT(\wTF)\bigr)_T
=
(-\ii\omega\vec a_0,\wTh)_\Omega
\qquad
\forall\wh\in\djdivz.
\end{multline}
The second equation enforces the discrete divergence-free constraint
\begin{equation}
\label{eq:expanded.div.constraint}
\sum_{T\in\t_h}
\bigl(
\DpT(\jTF),
q_{h\mid T}
\bigr)_T
=
0
\qquad
\forall q_h\in\Ph.
\end{equation}
\begin{remark}[Discrete current conservation]
\label{rmk:current.conservation}
By~\eqref{eq:Dph.zero.mean}, the constraint~\eqref{eq:expanded.div.constraint} implies $\Dph\jh=0$. Taking $q=1$ in~\eqref{eq:def.DT} on each cell therefore gives
\[
\sum_{F\in\f_T}\varepsilon_{T,F}\int_F\jFn=0
\qquad\forall T\in\t_h.
\]
Thus the shared face unknowns satisfy exact cellwise flux balance. This property concerns the hybrid current. The broken cell polynomial $\jTh$ is not required to have a continuous normal trace across faces.
\end{remark}

\subsection{Interpolation operators}
\label{sec:interpolation.operators}
We denote by $\pi_h^{\ell}:\LL\to\Poly^{\ell}(\t_h)$ the broken $L^2$-orthogonal projector defined cellwise by
\[
(\pi_h^{\ell}q)_{\mid T}
\defi
\pi_T^{\ell}(q_{\mid T})
\qquad
\forall T\in\t_h.
\]

We similarly denote by $\vproj{h}$ the cellwise vector-valued $L^2$ projector.
We define the reduction operators on fields with cellwise $H^1$ regularity, so their normal traces on individual faces belong to $L^2(F)$.

\noindent For each cell $T\in\t_h$, the local reduction operator is defined by
\begin{equation}
\label{eq:local.reduction}
\IhT:\vec H^1(T)\to\djdivT,
\qquad
\IhT(\vec v)
\defi
\left(
\vec\pi_T^{\ell}\vec v,
\bigl(
\pi_F^{\ell}(\vec v\cdot\normal_F)
\bigr)_{F\in\f_T}
\right).
\end{equation}
The face component uses the scalar $L^2(F)$ projection of the normal trace.

\noindent The global reduction operator is obtained by assembling the local cell and face components as
\begin{equation}
\label{eq:global.reduction}
\Ihh:\Hzdiv\cap\vec H^1(\t_h)\to\djdivz,
\qquad
\Ihh(\vec v)
\defi
\left(
\bigl(
\vec\pi_T^{\ell}(\vec v_{\mid T})
\bigr)_{T\in\t_h},
\bigl(
\pi_F^{\ell}(\vec v\cdot\normal_F)
\bigr)_{F\in\f_h}
\right).
\end{equation}
Here $\vec H^1(\t_h)$ denotes the broken Sobolev space. The homogeneous normal boundary condition ensures that $\Ihh\vec v\in\djdivz$ for every $\vec v\in\Hzdiv\cap\vec H^1(\t_h)$.

\noindent The reduction operator commutes with the divergence reconstruction in the sense that
\begin{equation}
\label{eq:global.commuting.Dph}
\Dph(\Ihh\vec v)
=
\pi_h^{\ell}(\Div\vec v)
\qquad
\forall\vec v\in\Hzdiv\cap\vec H^1(\t_h).
\end{equation}
Indeed, for every $T\in\t_h$ and $q\in\Poly^{\ell}(T)$, the definitions of $\DpT$ and $\IhT$, together with the orthogonality of the $L^2$ projectors and cellwise integration by parts, yield
\[
\bigl(\DpT(\IhT\vec v),q\bigr)_T
=
-(\vec v,\Grad q)_T
+
\sum_{F\in\f_T}
\varepsilon_{T,F}
(\vec v\cdot\normal_F,q_{\mid F})_F
=
(\Div\vec v,q)_T.
\]
Since $\DpT(\IhT\vec v)\in\Poly^\ell(T)$, this identity proves~\eqref{eq:global.commuting.Dph}.

\noindent Finally, the stabilization satisfies polynomial consistency. More precisely, for every $T\in\t_h$ and every $\vec p\in\bPoly^{\ell}(T)$, letting $\underline{\mathbbmsl p}_T\defi\IhT\vec p$, one has
\begin{equation}
\label{eq:poly.consistency.stab}
s_{\Div,T}(\underline{\mathbbmsl p}_T,\wTF)
=
0
\qquad
\forall\wTF\in\djdivT.
\end{equation}
In fact, the normal trace of $\vec p$ satisfies
$$
\vec p_{\mid F}\cdot\normal_F\in\Poly^\ell(F)
\qquad
\forall F\in\f_T,
$$
and therefore
$$\pi_F^\ell
\bigl( \vec p_{\mid F}\cdot\normal_F
\bigr) = \vec p_{\mid F}\cdot\normal_F,
$$
which implies \eqref{eq:poly.consistency.stab}.

\subsection{Norms and main results}
\label{sec:main.results}
The following results give the stability and approximation guarantees used to interpret the numerical study.
We equip $\djdivz$ with the discrete graph norm
\begin{equation}
\label{eq:def.discrete.norm}
\tnorm{h}{\jh}^2
\defi
\Norm{0,\Omega}{\jTh}^2
+
\Norm{0,\Omega}{\Dph\jh}^2
+
s_h(\jh,\jh),
\qquad
\tnorm{h}{(\jh,u_h)}^2
\defi
\tnorm{h}{\jh}^2
+
\Norm{0,\Omega}{u_h}^2.
\end{equation}

\begin{theorem}[Well-posedness and stability of the discrete problem]
\label{thm:well.posedness.discrete}
Assume the mesh-regularity conditions of Section~\ref{sec:discrete.spaces}, let $\mu>0$ be constant, let $\rho\in L^\infty(\Omega)$ be real-valued with $0<\rho_\ast\leq\rho\leq\rho^\ast$, and let $G$ be the Newtonian operator defined in Section~\ref{sec:continuous.model}. Then, for every fixed $\alpha>0$, the discrete mixed problem~\eqref{eq:discrete.problem} is well posed and admits a unique solution $(\jh,u_h)\in\djdivz\times\Ph$.
Moreover, there exists a constant $C>0$ such that
\begin{equation}
\label{eq:discrete.solution.stability}
\tnorm{h}{\jh}
+
\Norm{0,\Omega}{u_h}
\leq
C\,\omega\Norm{0,\Omega}{\vec a_0}.
\end{equation}
The constant may depend on the polynomial degree, the mesh-regularity parameters, the domain, the physical coefficients, $\omega$, and $\alpha$, but is independent of $h$.
\end{theorem}

\begin{theorem}[A priori error estimate]\label{thm:error.estimate}
Let $(\vec j,u)\in\Hzdiv\times\LLstar$ and $(\jh,u_h)\in\djdivz\times\Ph$ solve the continuous and discrete problems, respectively. Under the assumptions of Theorem~\ref{thm:well.posedness.discrete}, and assuming $u\in H^1(\Omega)$ and $\vec j|_T\in\vec H^1(T)$ for every cell,
\begin{align}\label{eq:error.estimate}
\tnorm{h}{(\jh,u_h)-(\Ihh\vec j,\proj{h}u)}
&\lesssim \Norm{0,\Omega}{\vec j-\vproj{h}\vec j}
 + \left(\sum_{T\in\t_h}\left[h_T^{-2}\Norm{0,T}{u-\proj{T}u}^{2}
       +\Norm{0,T}{\Grad(u-\proj{T}u)}^{2}\right]\right)^{1/2}\notag\\
&\quad + \left(\sum_{T\in\t_h}\sum_{F\in\f_T}h_F
           \Norm{0,F}{(\vec j-\vproj{T}\vec j)\cdot\normal_F}^{2}\right)^{1/2}.
\end{align}
If, in addition, $\vec j|_T\in\vec H^{r+1}(T)$ and $u|_T\in H^{r+1}(T)$ for every $T\in\t_h$ and some integer $0\leq r\leq\ell$, then
\begin{align}\label{eq:error.rate}
\tnorm{h}{(\jh,u_h)-(\Ihh\vec j,\proj{h}u)}
\lesssim\left(\sum_{T\in\t_h}h_T^{2(r+1)}|\vec j|_{\vec H^{r+1}(T)}^2\right)^{1/2} + \left(\sum_{T\in\t_h}h_T^{2r}|u|_{H^{r+1}(T)}^2\right)^{1/2}.
\end{align}
The hidden constants are independent of $h$, but may depend on the polynomial degree, the mesh-regularity parameters, the domain, the physical coefficients, and the stabilization parameter.
\end{theorem}

\noindent For sufficiently smooth solutions and fixed $\ell\geq1$, the current-projection and stabilization terms in~\eqref{eq:error.estimate} are of order $h^{\ell+1}$.
Estimating the trace of $u-\proj{T}u$ in the face residual gives order $h^\ell$, which determines the rate guaranteed in the discrete graph norm. By the triangle inequality and the approximation properties of $\vproj{h}$, the same $h^\ell$ bound holds for the cell current in $L^2(\Omega)^3$.
The proofs are given in Section~\ref{sec:stability.error}. The regularity assumptions on the current are local to each cell, so the stated rates apply across fitted material interfaces whenever the required cellwise regularity holds. The constants may depend on $\rho_\ast$ and $\rho^\ast$; no uniformity in the contrast $\rho^\ast/\rho_\ast$ is claimed.

\noindent For a prescribed coil excitation, the conductor contribution to the impedance is a bounded linear functional of the current. The following corollary therefore relates the discretization error in the impedance to the $L^2$ error in the cell current.

\begin{corollary}[Linear current observables]
\label{cor:current.observables}
Under the assumptions of Theorem~\ref{thm:error.estimate}, let $\vec b\in\bLL$ and define the complex-linear functional
\[
\mathcal Q_{\vec b}(\vec v)\defi\int_\Omega\vec b\cdot\vec v,
\]
where the dot product is bilinear, without complex conjugation. Then
\begin{align}
\bigl|\mathcal Q_{\vec b}(\vec j)-\mathcal Q_{\vec b}(\jTh)\bigr|
&\leq\Norm{0,\Omega}{\vec b}\Norm{0,\Omega}{\vec j-\jTh}\notag\\
&\leq\Norm{0,\Omega}{\vec b}\left(
\Norm{0,\Omega}{\vec j-\vproj{h}\vec j}
+\tnorm{h}{(\jh,u_h)-(\Ihh\vec j,\proj{h}u)}\right).
\label{eq:current.observable.error}
\end{align}
Consequently, Theorem~\ref{thm:error.estimate} gives an a priori bound for each such observable.
\end{corollary}
\begin{proof}
Apply the Cauchy--Schwarz inequality and use that the cell component of $\Ihh\vec j$ is $\vproj{h}\vec j$.
\end{proof}

\begin{remark}[Separated conducting components]
\label{rmk:separated.conductors}
The formulation and estimates extend to a fixed finite union $\Omega=\bigcup_{m=1}^M\Omega_m$ of pairwise separated, bounded, connected Lipschitz conducting components.
The scalar potential is determined up to one additive constant on each component. We therefore replace $\LLstar$ by
\[
  L^2_{\star,\mathrm c}(\Omega)
  \defi \left\{q\in L^2(\Omega)\;\middle|\;
    \int_{\Omega_m} q=0 \quad\forall m=1,\ldots,M \right\}
\]
and impose the same zero-mean conditions on the discrete multiplier space.
Writing $\xi_T$ for the integral of the prescribed divergence over $T$, the compatibility condition of Lemma~\ref{lem:relative.cochain.poincare} becomes $\sum_{T\subset\Omega_m}\xi_T=0$ for each $m$.
The discrete right inverse of Lemma~\ref{lem:stable.discrete.lifting} can then be constructed independently on each component. Summing the squared lifting bounds gives the global estimate.
The Green operator still acts on all of $\Omega$, so currents in distinct components remain magnetically coupled. Its boundedness and positivity hold on the full union, and the remaining stability and consistency arguments carry over with constants that may depend on the fixed component geometries.
The numerical experiments below use a connected conductor.
\end{remark}

\section{Implementation and ECT observables}
\label{sec:implementation}
The implementation is based on the HArDCore3D library\footnote{\url{https://github.com/jdroniou/HArDCore3D-release}}. We describe the algebraic structure of the mixed system, the quantities used to compare computational costs, and the quadrature used to assemble the nonlocal Green operator. We then define the ECT observations obtained from the computed current. The source code accompanying the numerical results provides the implementation of these steps.

\subsection{Algebraic structure and computational cost}
Only the cell-current unknowns enter the dense Green block. The face fluxes and the scalar multiplier enter through the local stabilization and divergence operators. Ordering the current coefficients by cell and then by face gives
\[
\begin{pmatrix}
\mathbf A_h & \mathbf B_h^*\\
\mathbf B_h & \mathbf 0
\end{pmatrix}
\begin{pmatrix}\mathbf J\\\mathbf U\end{pmatrix}
=
\begin{pmatrix}\mathbf f\\\mathbf 0\end{pmatrix},
\qquad
\mathbf A_h=\mathbf A_h^{\rm loc}+\ii\omega
\begin{pmatrix}\mathbf G_h&\mathbf 0\\\mathbf 0&\mathbf 0\end{pmatrix}.
\]
Here $\mathbf A_h^{\rm loc}$ contains the resistive and stabilization contributions, $\mathbf B_h$ represents the divergence coupling, and $\mathbf G_h$ represents $g_h$ on the cell-current space.
At fixed conductor geometry, mesh, polynomial degree, and permeability, the Green block is independent of the resistivity and the excitation frequency.
Changing the resistivity updates the local resistive contribution, while changing the frequency rescales the inductive contribution; these parameter changes generally require a new factorization of the full system.
Moving an air-cored source coil at fixed physical parameters changes only the excitation and observation functionals, so the full matrix and its factorization can be reused, as in the ECT benchmark.
A change in the geometry of an insulating defect generally requires updating the conductor mesh and recomputing the Green interactions involving the modified cells.

\noindent For $N_T$ cells and a uniform degree $\ell$, the number of cell-current unknowns entering the Green block is
\begin{equation}
\label{eq:def.NG}
N_G\defi\sum_{T\in\t_h}3\dim\Poly^\ell(T)
=3N_T\binom{\ell+3}{3}.
\end{equation}
With $N_F^{\rm int}$ interior faces, the number of independent unknowns after imposing the homogeneous boundary conditions and the zero-average multiplier constraint is
\begin{equation}
\label{eq:def.Ntot}
N_{\rm tot}=4N_T\binom{\ell+3}{3}
+N_F^{\rm int}\binom{\ell+2}{2}-1.
\end{equation}
For $M$ separated conducting components, the final term in~\eqref{eq:def.Ntot} is $-M$ rather than $-1$.
These counts measure different costs. $N_G$ is the dimension of the Green matrix block, while $N_{\rm tot}$ also includes the face and multiplier variables needed in the mixed solve.
Full storage of the Green block grows as $N_G^2$, up to symmetry and component storage conventions.
Increasing the polynomial degree enlarges this block on a fixed mesh.
Section~\ref{sec:dense.dof.efficiency} compares these mesh--degree choices using the cell-current error, Green block storage, and assembly and solution times.

\begin{remark}[Static condensation]
\label{rmk:static.condensation}
The Green term couples cell unknowns belonging to distinct mesh elements, so the cell--cell block is not block diagonal over the mesh. The standard HHO elimination of cell unknowns by independent elementwise solves is therefore unavailable. A global Schur complement remains possible, but requires operations on the globally coupled cell block. The hybrid construction retains its conservation and approximation properties without this usual local elimination advantage.

\end{remark}

\subsection{Quadrature of the Green operator}
\label{rmk:singular.self.interaction}
The assembly of the nonlocal form involves cell-pair contributions
\[
\int_T\int_{T'}
\frac{\boldsymbol\phi_T(\vec x)\cdot\boldsymbol\psi_{T'}(\vec y)}{4\pi|\vec x-\vec y|}
\,\mathrm d\vec y\,\mathrm d\vec x.
\]
The kernel is smooth for well-separated cells. Self-interactions contain an integrable singularity on the diagonal, while for distinct cells that touch the singularity is confined to coincident boundary points. Across a narrow gap the kernel remains smooth but can vary sharply.
For a discussion of singularity extraction with volume-uniform basis functions, see~\cite[Section~6.2]{Passarotto.Pitassi.Specogna:22}.

\noindent In the ECT implementation we use product quadrature for distinct cells. For self-interactions, we evaluate the inner integral by a radial transformation centred at each outer quadrature point.
Contributions from the cell faces are combined with their signed geometric factors.
The radial Jacobian removes the kernel singularity, and the resulting integrals are evaluated using face and radial quadrature.
Appendix~\ref{app:ect.setup} gives the oriented integration formula and the integration parameters, while Section~\ref{sec:numerics.fastener.hole} reports the Green quadrature sensitivity. The corresponding quadrature routines are provided in the source code accompanying the numerical results.

\noindent The analysis assumes exact integration; the effect of Green quadrature is assessed numerically in the ECT benchmark.

\subsection{ECT observation functionals}

\label{sec:ect.observables}
The computed current can be used to evaluate coil responses by reciprocity.
For a single coil carrying a current of amplitude $I>0$, let $\vec a_0$ be its prescribed excitation potential and $\jTh$ the resulting cell-current field.
With the phasor convention of~\eqref{eq:strong.model} and the voltage orientation adopted here, the conductor contribution to the coil impedance is
\begin{equation}\label{eq:reaction.impedance}
Z_{\mathrm{ec},h}=\frac{\ii\omega}{I^2}\int_\Omega
\vec a_0(\vec x)\cdot\jTh(\vec x)\,\mathrm d\vec x.
\end{equation}
The dot product is bilinear, without complex conjugation.
If $s$ denotes a probe position, Corollary~\ref{cor:current.observables} gives, for exact integration and prescribed excitation,
\begin{equation}
|Z_{\mathrm{ec}}(s)-Z_{\mathrm{ec},h}(s)|
\leq\frac{\omega}{I^2}\Norm{0,\Omega}{\vec a_0^{\,s}}
\Norm{0,\Omega}{\vec j^{\,s}-\jTh^{\,s}}.
\label{eq:impedance.error.bound}
\end{equation}
The same construction applies to distinct source and receiving coils.
Let $\jTh^{\mathrm s}$ be the response to a source coil carrying a current of amplitude $I_{\mathrm s}>0$, and let $\vec a_{\mathrm r}$ be the prescribed potential produced by a virtual current $I_{\mathrm r}>0$ in a receiving coil.
With the voltage orientation and phasor convention used in~\eqref{eq:reaction.impedance}, the conductor contribution to the transfer impedance is
\begin{equation}
Z_{\mathrm{rs},h}^{\mathrm{ec}}
\defi\frac{\ii\omega}{I_{\mathrm r}I_{\mathrm s}}
\int_\Omega\vec a_{\mathrm r}\cdot\jTh^{\mathrm s}.
\label{eq:transfer.impedance}
\end{equation}
Here receiving coils are treated as non-perturbing observations in the same magnetically homogeneous environment.
The single-coil response corresponds to coincident source and receiver, while subtracting two receiver functionals gives a differential signal.
These observations introduce no additional conductor unknowns and satisfy Corollary~\ref{cor:current.observables} whenever their prescribed potentials belong to $\bLL$.
The benchmark below uses the single-coil response.

\begin{remark}[Discrete dissipative response]
\label{rmk:discrete.dissipation}
Consider the exactly integrated scheme and a single coil with current amplitude $I>0$ and a real-valued excitation potential $\vec a_0$, as obtained by fixing the phase of a prescribed winding current.
The conductor impedance defined in~\eqref{eq:reaction.impedance} satisfies
\begin{equation}
I^2 Z_{\mathrm{ec},h}
= (\rho\jTh,\jTh)_\Omega
+\alpha s_h(\jh,\jh)
-\ii\omega(\mu G(\jTh),\jTh)_\Omega.
\label{eq:impedance.dissipation}
\end{equation}
Indeed, testing the first discrete equation with $\jh$ eliminates the multiplier term by the discrete divergence constraint and gives
$a_h(\jh,\jh)=-\ii\omega(\vec a_0,\jTh)_\Omega$.
Conjugating this equality proves~\eqref{eq:impedance.dissipation}.
In particular, $\Re Z_{\mathrm{ec},h}\geq0$; the stabilization contributes numerical dissipation in addition to the resistive term.
This sign property concerns the unsubtracted conductor impedance and imposes no sign on a reference-subtracted signal.
With quadrature, an exact algebraic identity requires the same assembled functional for excitation and readout and a Hermitian Green block; separate source and readout quadratures, as used below, can introduce a discrepancy.
\end{remark}

\section{Numerical results}
\label{sec:numerical.results}

The numerical experiments assess convergence, computational cost, and the prediction of an experimental ECT response.
A manufactured solution verifies convergence on several mesh families.
We then use the same problem to examine how refining the mesh or increasing the polynomial degree affects the current error, the size of the Green block, and the computational cost.
Finally, an experimental ECT benchmark assesses the impedance response on a locally resolved defective plate.

\noindent
For all tests with constant resistivity, the stabilization parameter is set to $\alpha=\rho$; thus $\alpha=1$ in the nondimensional manufactured problem and $\alpha=\sigma^{-1}$ in the physical benchmark.

\subsection{Manufactured solution on a cubic domain}
\label{sec:numerics.cube}
We first consider the nondimensional problem on $\Omega=(0,1)^3$ with
$\rho=\mu=\omega=\alpha=1$. The exact current and scalar multiplier are defined by
\begin{align*}
\psi(x,y,z)&=\sin^2(\pi x)\sin^2(\pi y)\sin^2(\pi z),\\
\boldsymbol j&=\operatorname{curl}(0,0,\psi)^{\mathsf T}
=(\partial_y\psi,-\partial_x\psi,0)^{\mathsf T},\\
u(x,y,z)&=\sin(2\pi x)\sin(2\pi y)\sin(2\pi z).
\end{align*}
By construction, $\operatorname{div}\boldsymbol j=0$ in $\Omega$,
$\boldsymbol j\cdot\boldsymbol n=0$ on $\partial\Omega$, and $\int_\Omega u=0$.
The excitation is prescribed as
\begin{equation}
\label{eq:manufactured.a0}
\vec a_0 = \frac{\ii}{\omega} \left( \rho\vec j+\Grad u \right) - \mu G(\vec j).
\end{equation}
With this choice, $(\vec j,u)$ solves the continuous problem exactly.
The prescribed excitation is evaluated using numerical integration of the Green potential in~\eqref{eq:manufactured.a0}.

\noindent We use Cartesian, tetrahedral, and Voronoi mesh families for the numerical simulations for polynomial degrees $\ell=1,2$. We report the relative errors
\begin{align*}
E_{j,0}&\defi\frac{\Norm{0,\Omega}{\vec j-\jTh}}{\Norm{0,\Omega}{\vec j}},\qquad
E_{u,0}\defi\frac{\Norm{0,\Omega}{u-u_h}} {\Norm{0,\Omega}{u}},\\
E_{j,h}&\defi \frac{ \tnorm{h}{\Ihh\vec j-\jh} }{ \tnorm{h}{\Ihh\vec j} }, \qquad r_{\Div,h} \defi \Norm{0,\Omega}{\Dph(\jh)}.
\end{align*}

\begin{figure}[!htbp]
\centering
\begin{tikzpicture}
\begin{groupplot}[
group style={group size=3 by 1,horizontal sep=1.1cm},
scale only axis,width=.255\textwidth,height=.22\textwidth,xmode=log,ymode=log,
xlabel={$N_T^{1/3}$},tick label style={font=\scriptsize},
label style={font=\small},title style={font=\small},
legend style={at={(.03,.03)},anchor=south west,font=\tiny,draw=none,fill=none},
every axis plot/.append style={thick,mark size=1.6pt}]
\nextgroupplot[xtick={2,4,8,16,32},xticklabels={2,4,8,16,32},xmin=1.76,xmax=35.84,ymin=0.00197021239,ymax=1.25028891,title={(a) Cartesian},ylabel={Relative error}]
\addplot[red,dashed,mark=*,mark options={fill=red}] table[x=resolution,y=Ej0,col sep=comma] {data/manufactured-cartesian-l1.csv};
\addlegendentry{$E_{j,0}$}
\addplot[blue,dashed,mark=square*,mark options={fill=blue}] table[x=resolution,y=Ejh,col sep=comma] {data/manufactured-cartesian-l1.csv};
\addlegendentry{$E_{j,h}$}
\addplot[green!60!black,dashed,mark=triangle*,mark options={fill=green!60!black}] table[x=resolution,y=Eu0,col sep=comma] {data/manufactured-cartesian-l1.csv};
\addlegendentry{$E_{u,0}$}

\draw[black,thin,solid] (axis cs:20.8,0.00256127611)--(axis cs:14.3448276,0.00256127611)--(axis cs:14.3448276,0.00371385035)--cycle;
\node[font=\scriptsize,anchor=east] at (axis cs:13.9144828,0.00371385035) {$1$};
\draw[black,thin,densely dotted] (axis cs:32.64,0.00256127611)--(axis cs:22.5103448,0.00256127611)--(axis cs:22.5103448,0.00538508301)--cycle;
\node[font=\scriptsize,anchor=east] at (axis cs:21.8350345,0.00538508301) {$2$};
\nextgroupplot[xtick={3,6,9,14},xticklabels={3,6,9,14},xmin=2.34819345,xmax=16.0174476,ymin=0.012278303,ymax=1.22444625,title={(b) Tetrahedral},yticklabel pos=left]
\addplot[red,dashed,mark=*,mark options={fill=red}] table[x=resolution,y=Ej0,col sep=comma] {data/manufactured-tetrahedral-l1.csv};
\addplot[blue,dashed,mark=square*,mark options={fill=blue}] table[x=resolution,y=Ejh,col sep=comma] {data/manufactured-tetrahedral-l1.csv};
\addplot[green!60!black,dashed,mark=triangle*,mark options={fill=green!60!black}] table[x=resolution,y=Eu0,col sep=comma] {data/manufactured-tetrahedral-l1.csv};

\draw[black,thin,solid] (axis cs:9.29584014,0.0159617939)--(axis cs:6.41092424,0.0159617939)--(axis cs:6.41092424,0.0231446012)--cycle;
\node[font=\scriptsize,anchor=east] at (axis cs:6.21859651,0.0231446012) {$1$};
\draw[black,thin,densely dotted] (axis cs:14.5873184,0.0159617939)--(axis cs:10.0602196,0.0159617939)--(axis cs:10.0602196,0.0335596717)--cycle;
\node[font=\scriptsize,anchor=east] at (axis cs:9.75841299,0.0335596717) {$2$};
\nextgroupplot[xtick={3,5,7,9,11},xticklabels={3,5,7,9,11},xmin=2.64,xmax=12.32,ymin=0.0142200869,ymax=1.41300753,title={(c) Voronoi},yticklabel pos=left]
\addplot[red,dashed,mark=*,mark options={fill=red}] table[x=resolution,y=Ej0,col sep=comma] {data/manufactured-voronoi-l1.csv};
\addplot[blue,dashed,mark=square*,mark options={fill=blue}] table[x=resolution,y=Ejh,col sep=comma] {data/manufactured-voronoi-l1.csv};
\addplot[green!60!black,dashed,mark=triangle*,mark options={fill=green!60!black}] table[x=resolution,y=Eu0,col sep=comma] {data/manufactured-voronoi-l1.csv};

\draw[black,thin,solid] (axis cs:7.15,0.018486113)--(axis cs:4.93103448,0.018486113)--(axis cs:4.93103448,0.0268048639)--cycle;
\node[font=\scriptsize,anchor=east] at (axis cs:4.78310345,0.0268048639) {$1$};
\draw[black,thin,densely dotted] (axis cs:11.22,0.018486113)--(axis cs:7.73793103,0.018486113)--(axis cs:7.73793103,0.0388670526)--cycle;
\node[font=\scriptsize,anchor=east] at (axis cs:7.5057931,0.0388670526) {$2$};
\end{groupplot}
\end{tikzpicture}
\caption{Computed relative errors for the manufactured solution at degree $\ell=1$ on Cartesian, tetrahedral, and Voronoi meshes. The solid and dotted triangles indicate reference orders $1$ and $2$, respectively.}
\label{fig:five-k1}
\end{figure}
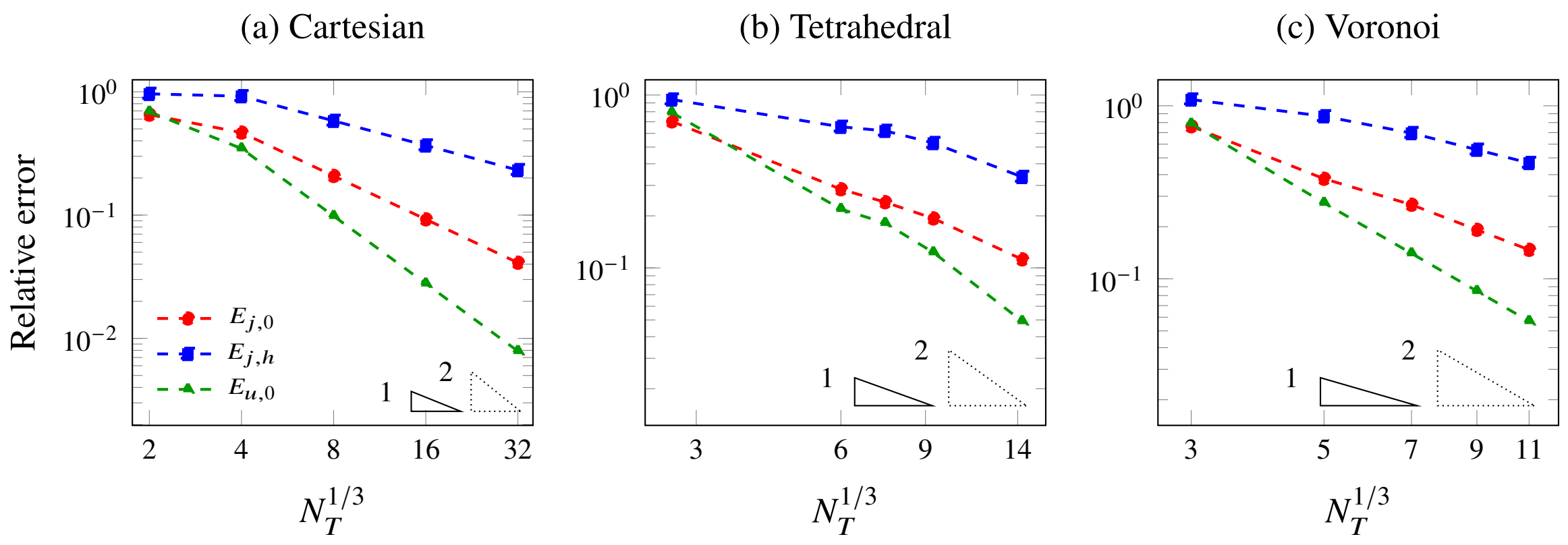

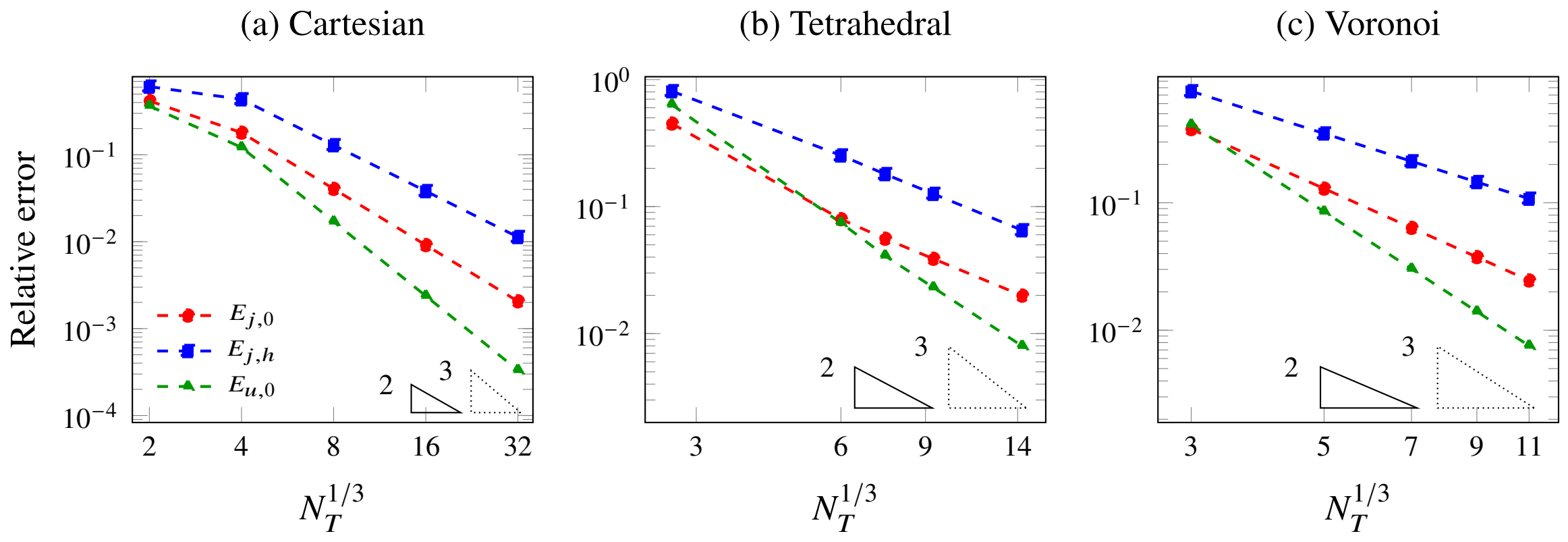
\begin{figure}[!htbp]
\centering
\begin{tikzpicture}
\begin{groupplot}[
group style={group size=3 by 1,horizontal sep=1.1cm},
scale only axis,width=.255\textwidth,height=.22\textwidth,xmode=log,ymode=log,
xlabel={$N_T^{1/3}$},tick label style={font=\scriptsize},
label style={font=\small},title style={font=\small},
legend style={at={(.03,.03)},anchor=south west,font=\tiny,draw=none,fill=none},
every axis plot/.append style={thick,mark size=1.6pt}]
\nextgroupplot[xtick={2,4,8,16,32},xticklabels={2,4,8,16,32},xmin=1.76,xmax=35.84,ymin=8.3018417e-05,ymax=0.791975423,title={(a) Cartesian},ylabel={Relative error}]
\addplot[red,dashed,mark=*,mark options={fill=red}] table[x=resolution,y=Ej0,col sep=comma] {data/manufactured-cartesian-l2.csv};
\addlegendentry{$E_{j,0}$}
\addplot[blue,dashed,mark=square*,mark options={fill=blue}] table[x=resolution,y=Ejh,col sep=comma] {data/manufactured-cartesian-l2.csv};
\addlegendentry{$E_{j,h}$}
\addplot[green!60!black,dashed,mark=triangle*,mark options={fill=green!60!black}] table[x=resolution,y=Eu0,col sep=comma] {data/manufactured-cartesian-l2.csv};
\addlegendentry{$E_{u,0}$}

\draw[black,thin,solid] (axis cs:20.8,0.000107923942)--(axis cs:14.3448276,0.000107923942)--(axis cs:14.3448276,0.000226910088)--cycle;
\node[font=\scriptsize,anchor=east] at (axis cs:13.9144828,0.000226910088) {$2$};
\draw[black,thin,densely dotted] (axis cs:32.64,0.000107923942)--(axis cs:22.5103448,0.000107923942)--(axis cs:22.5103448,0.000329019628)--cycle;
\node[font=\scriptsize,anchor=east] at (axis cs:21.8350345,0.000329019628) {$3$};
\nextgroupplot[xtick={3,6,9,14},xticklabels={3,6,9,14},xmin=2.34819345,xmax=16.0174476,ymin=0.00199600937,ymax=1.0532127,title={(b) Tetrahedral},yticklabel pos=left]
\addplot[red,dashed,mark=*,mark options={fill=red}] table[x=resolution,y=Ej0,col sep=comma] {data/manufactured-tetrahedral-l2.csv};
\addplot[blue,dashed,mark=square*,mark options={fill=blue}] table[x=resolution,y=Ejh,col sep=comma] {data/manufactured-tetrahedral-l2.csv};
\addplot[green!60!black,dashed,mark=triangle*,mark options={fill=green!60!black}] table[x=resolution,y=Eu0,col sep=comma] {data/manufactured-tetrahedral-l2.csv};

\draw[black,thin,solid] (axis cs:9.29584014,0.00259481219)--(axis cs:6.41092424,0.00259481219)--(axis cs:6.41092424,0.00545559262)--cycle;
\node[font=\scriptsize,anchor=east] at (axis cs:6.21859651,0.00545559262) {$2$};
\draw[black,thin,densely dotted] (axis cs:14.5873184,0.00259481219)--(axis cs:10.0602196,0.00259481219)--(axis cs:10.0602196,0.0079106093)--cycle;
\node[font=\scriptsize,anchor=east] at (axis cs:9.75841299,0.0079106093) {$3$};
\nextgroupplot[xtick={3,5,7,9,11},xticklabels={3,5,7,9,11},xmin=2.64,xmax=12.32,ymin=0.00188847691,ymax=0.975732017,title={(c) Voronoi},yticklabel pos=left]
\addplot[red,dashed,mark=*,mark options={fill=red}] table[x=resolution,y=Ej0,col sep=comma] {data/manufactured-voronoi-l2.csv};
\addplot[blue,dashed,mark=square*,mark options={fill=blue}] table[x=resolution,y=Ejh,col sep=comma] {data/manufactured-voronoi-l2.csv};
\addplot[green!60!black,dashed,mark=triangle*,mark options={fill=green!60!black}] table[x=resolution,y=Eu0,col sep=comma] {data/manufactured-voronoi-l2.csv};

\draw[black,thin,solid] (axis cs:7.15,0.00245501998)--(axis cs:4.93103448,0.00245501998)--(axis cs:4.93103448,0.00516167952)--cycle;
\node[font=\scriptsize,anchor=east] at (axis cs:4.78310345,0.00516167952) {$2$};
\draw[black,thin,densely dotted] (axis cs:11.22,0.00245501998)--(axis cs:7.73793103,0.00245501998)--(axis cs:7.73793103,0.0074844353)--cycle;
\node[font=\scriptsize,anchor=east] at (axis cs:7.5057931,0.0074844353) {$3$};
\end{groupplot}
\end{tikzpicture}
\caption{Computed relative errors for the manufactured solution at degree $\ell=2$ on Cartesian, tetrahedral, and Voronoi meshes. The solid and dotted triangles indicate reference orders $2$ and $3$, respectively.}
\label{fig:five-k2}
\end{figure}

\noindent Figures~\ref{fig:five-k1} and~\ref{fig:five-k2} report the computed errors at degrees $\ell=1$ and $2$. The errors decrease under mesh refinement, and the higher degree improves the cell-current and scalar-multiplier approximations. The reference slopes indicate the order $h^\ell$ guaranteed by Theorem~\ref{thm:error.estimate} and the order $h^{\ell+1}$. For the scalar multiplier, the observed $L^2$ convergence is close to $h^{\ell+1}$, one order higher than the rate guaranteed by the theorem.

\noindent The recorded reconstructed-divergence residuals range from $1.32\times10^{-15}$ to $1.94\times10^{-12}$, consistent with enforcement of the hybrid conservation constraint to linear-solve accuracy.

\FloatBarrier
\subsection{Effect of mesh size and polynomial degree on accuracy and cost}
\label{sec:dense.dof.efficiency}
We first compare finer degree-one meshes with coarser degree-two meshes to assess whether increasing the degree can reduce both the relative cell-current $L^2$ error $E_{j,0}$ and the size of the Green block.
We then compare discretizations from the same mesh family with similar Green block dimensions. The error is evaluated on the manufactured solution of Section~\ref{sec:numerics.cube}, and the degree is uniform within each run.

\noindent We distinguish the dimension $N_G$ of the Green block, defined in~\eqref{eq:def.NG}, from its storage $M_G$, the memory in GiB occupied by the assembled block, including its complex coefficients and matrix indices. In the uncompressed Eigen representation used here, zero component blocks are omitted and $M_G$ grows approximately as $N_G^2$. It does not include other parts of the mixed system or solver memory.
The time $t_G^{\rm ass}$ measures Green assembly; $t^{\rm sol}$ measures the first solution of the mixed system, including solver setup, LU factorization, the solution of the linear system, and residual checks. These first-solution costs are distinct from the reused-factorization solve times reported for the ECT scan in Section~\ref{sec:numerics.fastener.hole}.

\noindent Table~\ref{tab:vi-efficiency} exhibits several cases in which the coarser degree-two discretization gives a smaller error with fewer cell-current unknowns. On Cartesian meshes, $64$ degree-two cells give $E_{j,0}\simeq0.180$ with $N_G=1920$, compared with $E_{j,0}\simeq0.208$ and $N_G=6144$ for $512$ degree-one cells. The corresponding reported Green storage drops from $0.234$ to $0.023$~GiB. The selected tetrahedral and Voronoi pairs likewise show simultaneous reductions in error and nonlocal storage. Thus, on these smooth-solution tests, increasing the degree can reduce the Green block dimension required to meet an accuracy threshold.

\begin{table}[!htbp]
\centering
\caption{Accuracy and Green block storage. Each pair compares a finer mesh at $\ell=1$ (first row) with a coarser mesh at $\ell=2$ (second row). In the first three pairs, the coarser mesh gives a smaller error with fewer coupled current unknowns $N_G$; the last two pairs compare similar block dimensions within the same mesh family.}
\label{tab:vi-efficiency}
\small
\setlength{\tabcolsep}{6pt}
\begin{filecontents*}{data/mesh-degree-cost.csv}
family,degree,NT,NG,Ntot,Ej0,MG,assembly,solve
Cartesian,1,512,6144,12223,2.077e-1,0.234,1653.7,143.28
Cartesian,2,64,1920,3423,1.799e-1,0.023,1072.5,2.92
Tetrahedral,1,816,9792,17432,1.934e-1,0.595,79.1,424.05
Tetrahedral,2,216,6480,10847,7.963e-2,0.261,103.8,106.29
Voronoi,1,729,8748,25493,1.920e-1,0.475,65017.4,921.19
Voronoi,2,125,3750,8893,1.293e-1,0.087,38934.5,50.88
Tetrahedral,1,408,4896,8684,2.400e-1,0.149,33.4,51.81
Tetrahedral,2,216,6480,10847,7.963e-2,0.261,103.8,106.29
Voronoi,1,343,4116,11649,2.681e-1,0.105,23274.4,82.48
Voronoi,2,125,3750,8893,1.293e-1,0.087,38934.5,50.88
\end{filecontents*}
\pgfplotstabletypeset[
col sep=comma,
columns={family,degree,NT,NG,Ntot,Ej0,MG},
columns/family/.style={string type,column name={Family},column type=l},
columns/degree/.style={column name={$\ell$},column type=r},
columns/NT/.style={column name={Cells $N_T$},column type=r},
columns/NG/.style={column name={$N_G$},column type=r},
columns/Ntot/.style={column name={$N_{\rm tot}$},column type=r},
columns/Ej0/.style={column name={$E_{j,0}$},sci,sci zerofill,precision=3,column type=r},
columns/MG/.style={column name={$M_G$ (GiB)},fixed,fixed zerofill,precision=3,column type=r},
every head row/.style={before row=\toprule,after row=\midrule},
every nth row={2}{before row=\midrule},
every last row/.style={after row=\bottomrule},
1000 sep={}
]{data/mesh-degree-cost.csv}
\end{table}

\noindent The timing results in Table~\ref{tab:vi-efficiency-times} identify the limits of using $N_G$ alone as a cost predictor. Indeed, for the first tetrahedral pair, the degree-two Green assembly is slower despite the smaller interaction block, while the solution time is lower. In the last Voronoi pair, degree two improves accuracy and lowers storage, but again increases assembly time. Cell geometry, quadrature, and the mixed-system size all affect the total cost. The comparisons therefore measure the combined choice of mesh and degree in this implementation; their clearest benefit is the reduction of Green block dimension and storage at the reported accuracies.

\begin{table}[!htbp]
\centering
\caption{Assembly and first-solution times, in seconds, for the comparisons in Table~\ref{tab:vi-efficiency}, in the same order. The solution time includes solver setup, LU factorization, the solve, and residual checks.}
\label{tab:vi-efficiency-times}
\small
\setlength{\tabcolsep}{10pt}
\pgfplotstabletypeset[
col sep=comma,
columns={family,degree,NT,assembly,solve},
columns/family/.style={string type,column name={Family},column type=l},
columns/degree/.style={column name={$\ell$},column type=r},
columns/NT/.style={column name={Cells $N_T$},column type=r},
columns/assembly/.style={column name={$t_G^{\rm ass}$},fixed,fixed zerofill,precision=1,column type=r},
columns/solve/.style={column name={$t^{\rm sol}$},fixed,fixed zerofill,precision=2,column type=r},
every head row/.style={before row=\toprule,after row=\midrule},
every nth row={2}{before row=\midrule},
every last row/.style={after row=\bottomrule},
1000 sep={}
]{data/mesh-degree-cost.csv}
\end{table}

\begin{figure}[!htbp]
\centering
\begin{tikzpicture}
\begin{axis}[
width=.85\textwidth,height=.48\textwidth,xmode=log,ymode=log,
xlabel={Green block dimension $N_G$},ylabel={$E_{j,0}$},
legend columns=2,legend style={at={(.5,1.04)},anchor=south,font=\scriptsize,draw=none},
every axis plot/.append style={thick,mark size=2pt}]
\addplot[blue,dashed,mark=*] table[x=NG,y=Ej0,col sep=comma] {data/manufactured-cartesian-l1.csv};
\addlegendentry{Cartesian, $\ell=1$}
\addplot[red,dashed,mark=*] table[x=NG,y=Ej0,col sep=comma] {data/manufactured-cartesian-l2.csv};
\addlegendentry{Cartesian, $\ell=2$}
\addplot[blue,dashed,mark=triangle*] table[x=NG,y=Ej0,col sep=comma] {data/manufactured-tetrahedral-l1.csv};
\addlegendentry{Tetrahedral, $\ell=1$}
\addplot[red,dashed,mark=triangle*] table[x=NG,y=Ej0,col sep=comma] {data/manufactured-tetrahedral-l2.csv};
\addlegendentry{Tetrahedral, $\ell=2$}
\addplot[blue,dashed,mark=square*] table[x=NG,y=Ej0,col sep=comma] {data/manufactured-voronoi-l1.csv};
\addlegendentry{Voronoi, $\ell=1$}
\addplot[red,dashed,mark=square*] table[x=NG,y=Ej0,col sep=comma] {data/manufactured-voronoi-l2.csv};
\addlegendentry{Voronoi, $\ell=2$}
\end{axis}
\end{tikzpicture}
\caption{Computed cell-current error $E_{j,0}$ versus the dimension $N_G$ of the nonlocal interaction space. Colour distinguishes the degree and markers distinguish the mesh family.}
\label{fig:vi-efficiency}
\end{figure}
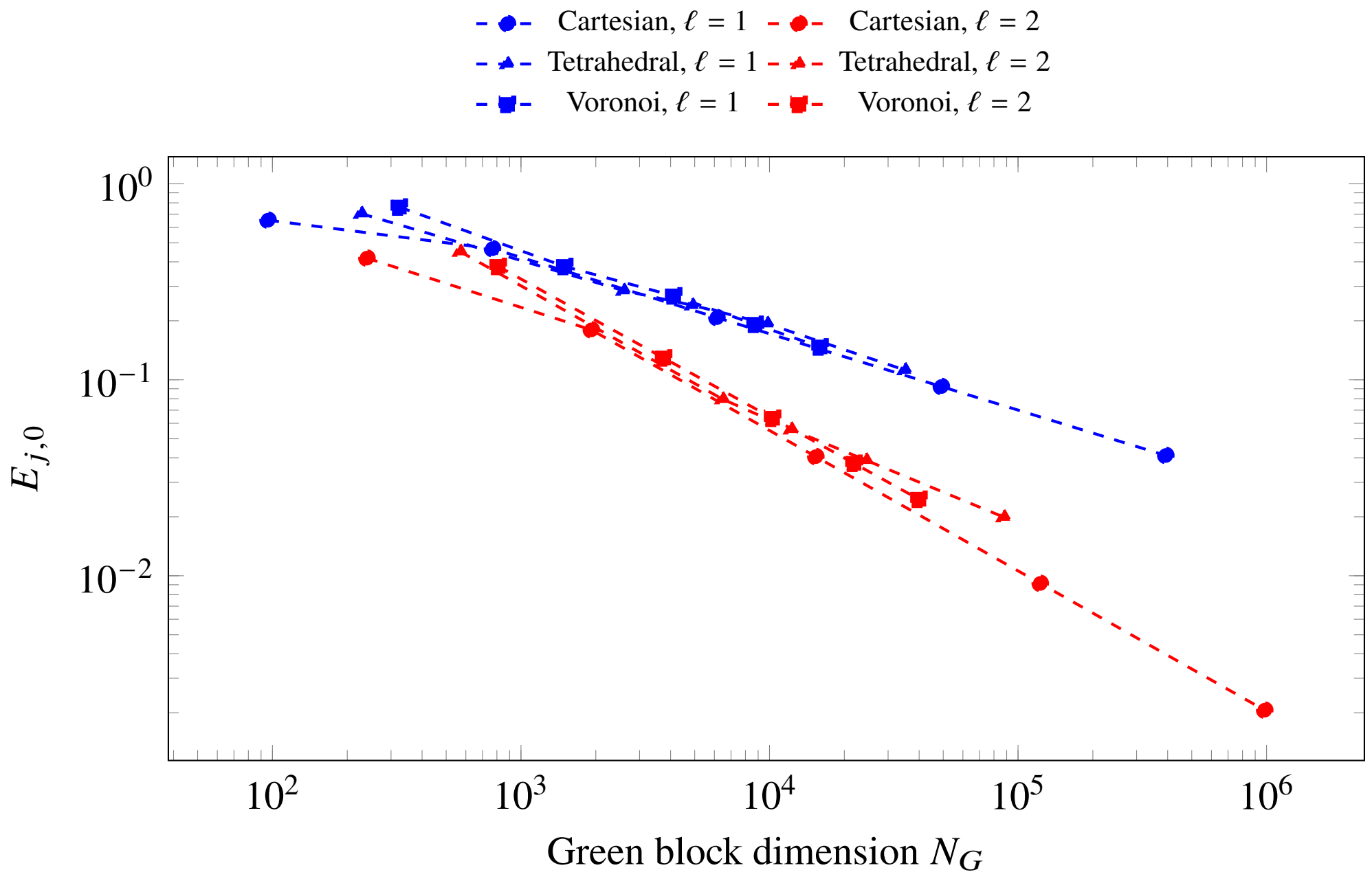

\noindent Figure~\ref{fig:vi-efficiency} places these comparisons within the available refinement sequences. In particular, the Voronoi degree-two point with $N_G=3750$ has error $0.129$, compared with $0.268$ at $N_G=4116$ for degree one. For this smooth solution, increasing the polynomial degree on a coarser mesh can thus reduce the cell-current error at a comparable Green block dimension. The degree is uniform within each run, and its choice depends on the target accuracy and the mesh family.

\FloatBarrier
\subsection{Eddy-current testing of a defective plate}
\label{sec:numerics.fastener.hole}
This example examines the impedance response of a defective conducting plate under a moving probe.
We first introduce the geometry and measured quantity, then compare locally refined polyhedral discretizations with the experimental scan and examine the resulting current field and computational costs.

\subsubsection{Benchmark and measured quantity}
We consider the single-plate experiment of~\cite[Test Case~1, Table~1.4]{Martinos.Theodoulidis.ea:13}, corresponding to configuration~B in~\cite{Martinos.Theodoulidis.ea:14}. An air-cored coil is scanned above an aluminium plate containing a circular fastener hole and a narrow radial through-thickness slot. The hole and the slot perturb the induced current paths, changing the reaction field and hence the coil impedance. This benchmark, also studied numerically in~\cite{Barbato.Poulakis.ea:15}, provides an experimental assessment of the method on a geometry with strongly separated length scales.

\noindent Table~\ref{tab:ect.parameters} collects the main parameters. The hole is centred at $(0,0)$, the plate occupies $[-150,150]^2\times[-2,0]$ in millimetres before removal of the hole and slot, and the slot extends along the positive $x$-axis to $x=19.8$~mm. The normal current vanishes on the entire conductor boundary. At $5$~kHz the skin depth is approximately $1.71$~mm, comparable to the plate thickness. All data are converted to SI units before assembly. Only the conductor is meshed; the coil enters through its prescribed excitation potential, detailed in Appendix~\ref{app:ect.setup}. The hole and slot are excluded from the conducting domain, so their insulating character is enforced geometrically. The common permeability $\mu_0$ in the aluminium and surrounding air places this configuration within the magnetically homogeneous model.

\begin{table}[!htbp]
\centering
\caption{Parameters of the single-plate ECT benchmark. Geometric dimensions are in millimetres.}
\label{tab:ect.parameters}
\small
\begin{tabular}{ll}
\hline
Plate and defect & $300\times300\times2$; hole radius $10$; slot $9.8\times0.234\times2$ \\
Material & $\sigma=17.34~\mathrm{MS/m}$, $\rho=\sigma^{-1}$, $\mu=\mu_0$ \\
Coil M1650 & $1650$ turns; inner/outer radii $7/12$; height $4$ \\
Excitation & $f=5~\mathrm{kHz}$, $I=1~\mathrm A$, lift-off $1.082$ \\
\hline
\end{tabular}
\end{table}
\noindent The $76$ experimental scan positions are $s_m=2m$~mm, $m=0,\ldots,75$, with coil centres $(x_{\rm c},y_{\rm c})=(s_m-70.8,0)$~mm. Thus, the hole and slot centres correspond to $s=70.8$ and $s=85.7$~mm, respectively. At each position $s_m$, the conductor contribution $Z_{\mathrm{ec},h}(s_m)$ is evaluated using~\eqref{eq:reaction.impedance} with the corresponding excitation potential and computed current.
Since the measurements contain the total coil impedance, we subtract the first scan position from both numerical and experimental values to obtain
\[
\Delta Z_h(s_m)=Z_{{\rm ec},h}(s_m)-Z_{{\rm ec},h}(s_0),\qquad
\Delta Z_{\rm exp}(s_m)=Z_{\rm exp}(s_m)-Z_{\rm exp}(s_0).
\]
The curves are normalized by the common reactance $X_0=\omega L_0\simeq1685.622~\Omega$, where $L_0=53.655$~mH is the measured free-space inductance. This subtraction removes the position-independent coil contribution and retains the impedance variation along the scan. The absolute discretization error in the reference-subtracted signal is bounded by the sum of the right-hand sides of~\eqref{eq:impedance.error.bound} at $s_m$ and $s_0$. These a priori bounds concern discretization of the stated model; the experimental discrepancy also reflects modelling and measurement effects.

\subsubsection{Polyhedral discretizations and impedance response}
We use three meshes, denoted by L0, L1, and L2, with $295$, $942$, and $2838$ polyhedral cells, respectively. The meshes are obtained by extruding polygonal partitions through one, two, and three thickness layers. Local resolution is concentrated around the hole and slot, while larger cells cover the remainder of the plate; Figure~\ref{fig:ect.mesh} illustrates this distribution. All levels retain the same polygonal approximation of the hole and the prescribed slot width. Additional mesh data are given in Appendix~\ref{app:ect.setup}. Table~\ref{tab:ect.results} reports the Green block dimension $N_G$ and the total number of independent unknowns $N_{\rm tot}$.

\paragraph{Scan setup.}
The L0 runs used an Intel Core Ultra 7 165U workstation with 15~GiB of RAM. L1 and L2 ran on Montpellier MESO nodes with AMD EPYC 9654 processors (192 physical cores and approximately 1.4~TiB of RAM per node). Each cluster run was allocated one CPU, with 128~GB of memory requested for L1 and 512~GB for L2. The HArDCore3D commit recorded for the cluster runs was \texttt{73fc1e9300a1}.

\noindent All three production scans used polynomial degree $\ell=1$ and Green quadrature parameter $q_G=4$. Distinct-cell interactions used products of degree-four cell rules; self-interactions used degree-eight cell and face rules with five-point Gauss--Legendre radial quadrature. Coil-source integration used a tensor-product midpoint rule with $4\times32\times4$ nodes in the radial, azimuthal, and axial directions. The cell quadrature degrees were eight for the source projection and ten for the impedance functional. Local polynomial terms were evaluated using exact monomial integration.

\input{figures/ect-mesh.tex}

\begin{figure}[!htbp]
\centering
\begin{tikzpicture}
\begin{groupplot}[group style={group size=2 by 1,horizontal sep=1.8cm},
scale only axis,width=.37\textwidth,height=.28\textwidth,scaled y ticks=false,yticklabel style={/pgf/number format/fixed,/pgf/number format/precision=2},xmin=0,xmax=150,xlabel={$s$ (mm)},
tick label style={font=\scriptsize},label style={font=\small},title style={font=\small},
every axis plot/.append style={line width=.8pt}]
\nextgroupplot[ytick={-0.02,-0.01,0,0.01,0.02},ylabel={$\Re\Delta Z/X_0$},title={Real part},
legend to name=ectscanlegend,legend columns=4,
legend style={font=\scriptsize,draw=none,column sep=.65em}]
\addplot[blue!65,densely dotted,no marks] table[x=s_mm,y=L0_re,col sep=comma] {data/ect-scan.csv};
\addlegendentry{L0}
\addplot[orange!90!black,dashed,no marks] table[x=s_mm,y=L1_re,col sep=comma] {data/ect-scan.csv};
\addlegendentry{L1}
\addplot[red!75!black,solid,no marks] table[x=s_mm,y=L2_re,col sep=comma] {data/ect-scan.csv};
\addlegendentry{L2}
\addplot[black,only marks,mark=o,mark size=1.1pt] table[x=s_mm,y=experiment_re,col sep=comma] {data/ect-scan.csv};
\addlegendentry{Experiment}
\nextgroupplot[ytick={0,0.05,0.1,0.15},ylabel={$\Im\Delta Z/X_0$},title={Imaginary part}]
\addplot[blue!65,densely dotted,no marks] table[x=s_mm,y=L0_im,col sep=comma] {data/ect-scan.csv};
\addplot[orange!90!black,dashed,no marks] table[x=s_mm,y=L1_im,col sep=comma] {data/ect-scan.csv};
\addplot[red!75!black,solid,no marks] table[x=s_mm,y=L2_im,col sep=comma] {data/ect-scan.csv};
\addplot[black,only marks,mark=o,mark size=1.1pt] table[x=s_mm,y=experiment_im,col sep=comma] {data/ect-scan.csv};
\end{groupplot}
\end{tikzpicture}\par\smallskip
\pgfplotslegendfromname{ectscanlegend}\par
\caption{Reference-subtracted impedance scans at $5$~kHz and degree $\ell=1$, normalized by the common free-space reactance. Experimental values are from Table~1.4 of~\cite{Martinos.Theodoulidis.ea:13}.}
\label{fig:ect.scan}
\end{figure}
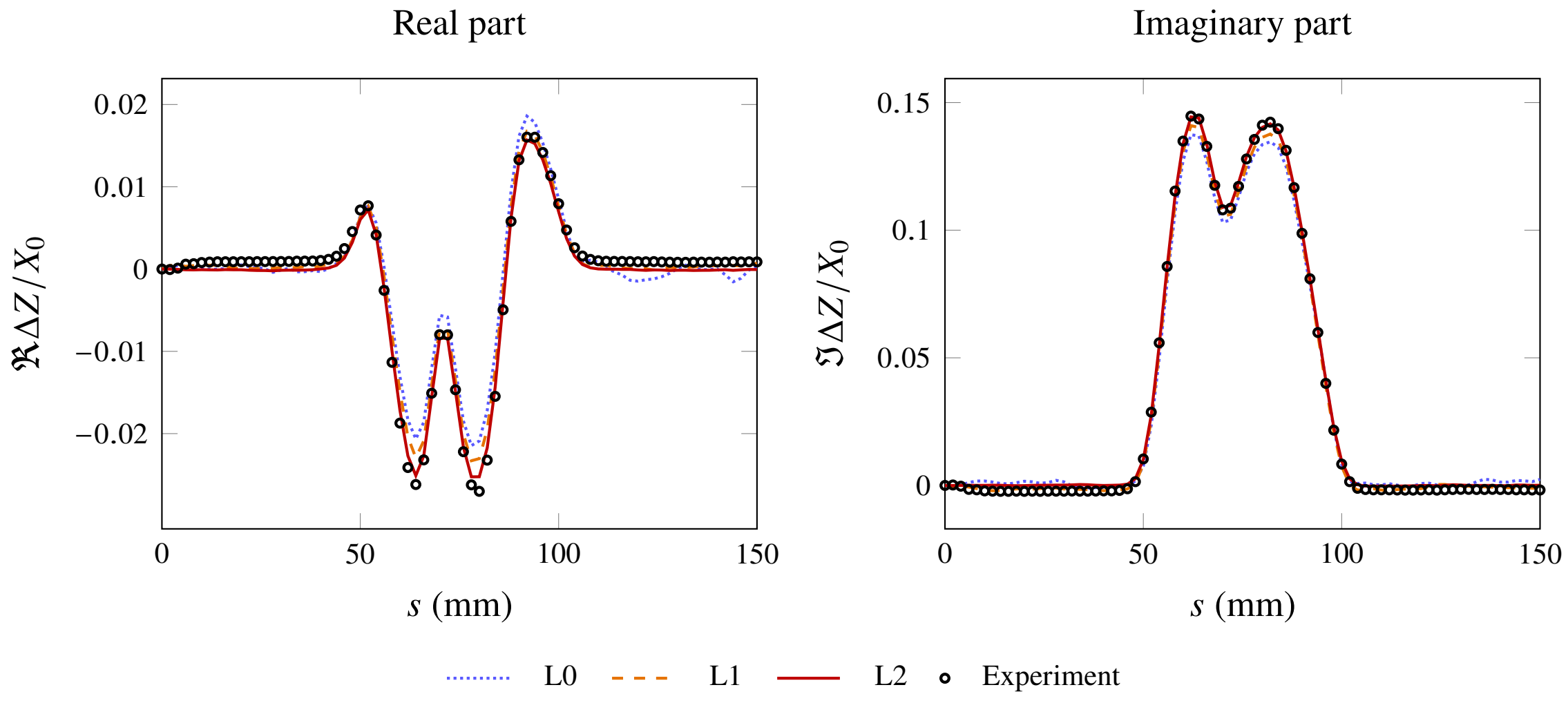

\paragraph{Impedance scans and numerical sensitivity.}
Figure~\ref{fig:ect.scan} compares the numerical and experimental real and imaginary parts of $\Delta Z/X_0$. To quantify the agreement, let $\boldsymbol{\Delta Z}_h=(\Delta Z_h(s_m))_{m=0}^{75}$ and define
\[
\varepsilon_2=\frac{\|\boldsymbol{\Delta Z}_h-\boldsymbol{\Delta Z}_{\rm exp}\|_{\ell^2}}{\|\boldsymbol{\Delta Z}_{\rm exp}\|_{\ell^2}},\qquad
\eta_r=\frac{\|\boldsymbol{\Delta Z}_{h_r}-\boldsymbol{\Delta Z}_{h_{r-1}}\|_{\ell^2}}{\|\boldsymbol{\Delta Z}_{h_r}\|_{\ell^2}}.
\]
The componentwise discrepancies $\varepsilon_{\Re}$ and $\varepsilon_{\Im}$ are defined as $\varepsilon_2$ with the real or imaginary part applied to both the numerator and denominator vectors. All reported quantities use the same reference subtraction and normalization.

\begin{table}[!htbp]
\centering
\caption{ECT discretizations and relative scan discrepancies at $\ell=1$. The last four columns are percentages; $\eta_r$ compares a mesh with the preceding level.}
\label{tab:ect.results}
\small
\setlength{\tabcolsep}{5pt}
\begin{filecontents*}{data/ect-discrepancies.csv}
mesh,NT,NG,Ntot,error_complex_percent,error_real_percent,error_imag_percent,mesh_difference_percent
L0,295,3540,7251,7.11,25.14,6.10,---
L1,942,11304,24776,3.74,14.65,3.06,4.05
L2,2838,34056,76841,2.94,10.05,2.55,3.67
\end{filecontents*}
\pgfplotstabletypeset[
col sep=comma,
columns/mesh/.style={string type,column name={Mesh},column type=l},
columns/NT/.style={string type,column name={$N_T$},column type=r},
columns/NG/.style={string type,column name={$N_G$},column type=r},
columns/Ntot/.style={string type,column name={$N_{\rm tot}$},column type=r},
columns/error_complex_percent/.style={string type,column name={$\varepsilon_2$},column type=r},
columns/error_real_percent/.style={string type,column name={$\varepsilon_{\Re}$},column type=r},
columns/error_imag_percent/.style={string type,column name={$\varepsilon_{\Im}$},column type=r},
columns/mesh_difference_percent/.style={string type,column name={$\eta_r$},column type=r},
every head row/.style={before row=\toprule,after row=\midrule},
every last row/.style={after row=\bottomrule}
]{data/ect-discrepancies.csv}
\end{table}

\noindent The discrepancy with the measurements decreases under refinement in both components, with a complex relative discrepancy of $2.94\%$ on L2. The corresponding maximum absolute discrepancy, divided by the maximum experimental signal magnitude, is $1.93\%$. The imaginary part dominates the complex signal; reporting the components separately shows the larger remaining relative discrepancy in the real part. At $s=86$~mm, with the coil centred at $(15.2,0)$~mm near the slot centre, L2 gives $\Delta Z_h=-6.84+220.69\ii~\Omega$, compared with $-8.34+221.30\ii~\Omega$ experimentally. The successive-mesh differences of $4.05\%$ and $3.67\%$ indicate that a mesh dependence remains at the reported resolutions.

\noindent Increasing $q_G$ from $4$ to $5$ changes the full L0 scan by $0.061\%$ in relative $\ell^2$ norm. On L1 and L2, the same check is performed at the reference position and at $s=86$~mm; the resulting reference-subtracted impedances change by $0.024\%$ and $0.012\%$, respectively. These changes are smaller than the measured mesh differences. The coil-source quadrature is held fixed in these checks, which therefore assess sensitivity to the Green integration.

\subsubsection{Current distribution and computational cost}
Figure~\ref{fig:ect.current} shows the real part of the L2 cell-current field near the slot for the coil centred at $(15.2,0)$~mm. The real current follows the slot sides and turns around its tip, where the sampled magnitudes are largest. Together with the impedance scan, this illustrates the use of local geometric resolution within a conductor-only volume integral discretization. The effect of increasing the polynomial degree is assessed separately in Section~\ref{sec:dense.dof.efficiency}.

\input{figures/ect-current.tex}

\noindent At fixed frequency, the system matrix is independent of the coil position. We therefore assemble and factorize it once, then update the excitation and evaluate the impedance at each position. Table~\ref{tab:ect.cost} separates assembly and factorization from the average solve time for the remaining $75$ right-hand sides. The direct solver is Eigen SparseLU, with iterative refinement. The implementation enforces the zero-average condition using one additional Lagrange multiplier, so the augmented systems have $N_{\rm tot}+2$ rows.

\begin{table}[!htbp]
\centering
\caption{Recorded ECT costs. Times are in seconds and peak resident memory is in GiB. The mean solve time excludes the first solve and reuses its factorization; it excludes excitation assembly and impedance evaluation. L0 was run locally, whereas L1 and L2 used single-CPU cluster allocations, so the timings are not a controlled hardware-scaling comparison.}
\label{tab:ect.cost}
\small
\begin{filecontents*}{data/ect-cost.csv}
mesh,assembly_seconds,factorization_seconds,mean_reused_solve_seconds,peak_memory_GiB
L0,670.9,38.8,0.307,0.942
L1,3110.5,1063.4,7.069,8.764
L2,12288.5,50170.0,85.012,86.346
\end{filecontents*}
\pgfplotstabletypeset[
col sep=comma,
columns/mesh/.style={string type,column name={Mesh},column type=l},
columns/assembly_seconds/.style={string type,column name={Assembly},column type=r},
columns/factorization_seconds/.style={string type,column name={Factorization},column type=r},
columns/mean_reused_solve_seconds/.style={string type,column name={Mean reused solve},column type=r},
columns/peak_memory_GiB/.style={string type,column name={Peak memory},column type=r},
every head row/.style={before row=\toprule,after row=\midrule},
every last row/.style={after row=\bottomrule}
]{data/ect-cost.csv}
\end{table}
\noindent On L2, matrix factorization took 13.9 hours and was performed once for the entire 76-position scan. The subsequent 75 linear solves took 85 seconds each on average, excluding excitation assembly and impedance evaluation.

\noindent The polyhedral mesh resolves the hole and narrow slot locally while retaining larger cells elsewhere. Since the volume integral formulation discretizes only the conductor, the probe can move over the plate without remeshing. The refinement around the slot stays fixed throughout the scan.
The benchmark tests the geometric representation of an insulating defect. Discontinuous resistivities are covered by the formulation and analysis but are not varied in these experiments.
The source code, meshes and input files used for the ECT benchmark are available from the authors upon request. The experimental measurements are reported in Table~1.4 of~\cite{Martinos.Theodoulidis.ea:13}.

\FloatBarrier
\section{Stability and error analysis}
\label{sec:stability.error}
\label{sec:stability.error.analysis}
We prove the results stated in Section~\ref{sec:main.results}.
A bounded right inverse of the reconstructed divergence supplies the discrete inf-sup condition, while the resistive and stabilization terms give coercivity on its kernel.
We then estimate the residual obtained by inserting the interpolated exact solution into the discrete equations and combine this consistency bound with mixed stability.

\subsection{Discrete graph norm and continuity}
\label{sec:discrete.graph.norm}
We first introduce the hybrid $L^2$-norm on $\djdivz$, defined by
\begin{equation}
\label{eq:hybrid.L2.norm}
\Norm{0,h}{\jh}^2
\defi
\sum_{T\in\t_h}
\left(
\Norm{0,T}{\jT}^2
+
\sum_{F\in\f_T}
h_F\Norm{0,F}{\jFn}^2
\right)
\qquad
\forall\jh\in\djdivz.
\end{equation}
Polynomial inverse trace inequalities and the mesh-regularity assumptions yield the uniform norm equivalence
\begin{equation}
\label{eq:hybrid.L2.stab.equiv}
\Norm{0,h}{\jh}^2
\simeq
\Norm{0,\Omega}{\jTh}^2
+
s_h(\jh,\jh)
\qquad
\forall\jh\in\djdivz,
\end{equation}
where the constants hidden in the notation $\simeq$ are independent of the mesh size.

\noindent It follows from~\eqref{eq:hybrid.L2.stab.equiv} that
\begin{equation}
\label{eq:graph.hybrid.equiv}
\tnorm{h}{\jh}^2
\simeq
\Norm{0,h}{\jh}^2
+
\Norm{0,\Omega}{\Dph\jh}^2.
\end{equation}
Thus, $\tnorm{h}{\cdot}$ is a discrete counterpart of the $\vec H(\Div)$ graph norm, whose $L^2$ component is represented by the hybrid norm~\eqref{eq:hybrid.L2.norm}.

\noindent Under the coefficient assumptions of Section~\ref{sec:continuous.model}, the properties of $G$ stated in Section~\ref{sec:continuous.model} yield constants $C_a,C_b>0$, independent of $h$, such that
\begin{align}
|a_h(\jh,\wh)|
&\leq
C_a\,
\tnorm{h}{\jh}\,
\tnorm{h}{\wh}
&&
\forall\jh,\wh\in\djdivz,
\label{eq:continuity.ah}
\\
|b_h(\wh,q_h)|
&\leq
C_b\,
\tnorm{h}{\wh}\,
\Norm{0,\Omega}{q_h}
&&
\forall\wh\in\djdivz,
\quad
\forall q_h\in\Ph.
\label{eq:continuity.bh}
\end{align}
The constants may depend on the polynomial degree, the mesh-regularity parameters, $\rho^\ast$, $\mu$, $\omega$, $\alpha$, and the operator norm of $G$, but they are independent of the mesh size.

\subsection{A relative top-degree cochain lifting}
\label{sec:relative.cochain.lifting}
We state the relative top-degree counterpart of the cochain Poincar\'e inequality proved in~\cite[Lemma~6]{DiPietro.Droniou.Hanot.Pitassi:25}.
It is stated directly in the three-dimensional setting required in the subsequent analysis.

\begin{lemma}[Relative top-degree cochain Poincar\'e inequality]
\label{lem:relative.cochain.poincare}
Let $(\xi_T)_{T\in\t_h}$ be a family of complex numbers satisfying
\begin{equation}
\sum_{T\in\t_h}\xi_T=0. \label{eq:relative.cochain.compatibility}
\end{equation}
Then there exists a family of complex numbers $(\lambda_F)_{F\in\f_h}$, depending linearly on $(\xi_T)_{T\in\t_h}$, such that
\begin{align}
\lambda_F&=0
\qquad
\forall F\in\f_h\cap\partial\Omega,\label{eq:relative.cochain.boundary}\\
\sum_{F\in\f_T}
\varepsilon_{T,F}\lambda_F
&=
\xi_T
\qquad
\forall T\in\t_h, \label{eq:relative.cochain.identity}
 \end{align}
 and
 \begin{equation}
\sum_{F\in\f_h}
h_F^{-1}|\lambda_F|^2
\lesssim
\sum_{T\in\t_h}
h_T^{-3}|\xi_T|^2, \label{eq:relative.cochain.bound}
\end{equation}
where the hidden constant is independent of $h$.
\end{lemma}
\noindent The proof is given in Appendix~\ref{app:cochain.proof}.
It transfers a stable Raviart--Thomas lifting on an auxiliary simplicial mesh to polyhedral face fluxes. The result therefore follows from established mixed finite element stability, while its relative cochain form incorporates the zero boundary fluxes and the associated compatibility condition. The bound inherits the finite element lifting and mesh-regularity constants. For separated conductors the lemma applies on each component under the compatibility condition stated in Remark~\ref{rmk:separated.conductors}.

\noindent To construct a uniformly bounded right inverse of the reconstructed divergence, we first use Lemma~\ref{lem:relative.cochain.poincare} to lift the cell integrals of the prescribed divergence to face fluxes. This step follows the top-degree construction on the mesh faces in~\cite{Pitassi:26}.
We then determine the cell unknowns by a local $L^2(T)$-Riesz representation on $\Grad\Poly^\ell(T)$.

\begin{lemma}[Stable discrete lifting of the reconstructed divergence]
\label{lem:stable.discrete.lifting} For every $q_h\in\Ph$, there exists $\underline{\mathbbmsl w}_{q_h}\in\djdivz$ and a constant $C>0$, independent of $h$, such that
\[
\Dph(\underline{\mathbbmsl w}_{q_h})
=
q_h,
\qquad
\tnorm{h}{\underline{\mathbbmsl w}_{q_h}}
\leq
C\Norm{0,\Omega}{q_h}.
\]
The lifting can be selected linearly in $q_h$.
\end{lemma}

\begin{proof}
Let $q_h\in\Ph$ and set
\[
\xi_T
\defi
\int_Tq_h
\qquad
\forall T\in\t_h.
\]
Since $(q_h,1)_\Omega=0$, we have
\[
\sum_{T\in\t_h}\xi_T=0.
\]
Applying Lemma~\ref{lem:relative.cochain.poincare} to $(\xi_T)_{T\in\t_h}$ yields a family $(\lambda_F)_{F\in\f_h}\subset\mathbb C$, depending linearly on $q_h$ and satisfying~\eqref{eq:relative.cochain.boundary}--\eqref{eq:relative.cochain.bound}.
Moreover, Cauchy--Schwarz, $|T|\simeq h_T^3$, and~\eqref{eq:relative.cochain.bound} give
\begin{equation}
\label{eq:cochain.lifting.bound}
\sum_{F\in\f_h}
h_F^{-1}|\lambda_F|^2
\lesssim
\Norm{0,\Omega}{q_h}^2.
\end{equation}

\noindent For every face $F\in\f_h$, define the constant polynomial
\begin{equation}
\label{eq:face.lifting.definition}
w_{F,n}
\defi
\frac{\lambda_F}{|F|}
\in
\Poly^{\ell}(F).
\end{equation}
In particular, $w_{F,n}=0$ on every boundary face. Moreover, using $|F|\simeq h_F^2$ in~\eqref{eq:cochain.lifting.bound}, we obtain
\begin{equation}
\label{eq:face.lifting.bound}
\sum_{F\in\f_h}
h_F\Norm{0,F}{w_{F,n}}^2
\lesssim
\Norm{0,\Omega}{q_h}^2.
\end{equation}
\noindent We now complete the face unknowns by constructing a cell polynomial. For each $T\in\t_h$, consider the functional on $\Grad\Poly^{\ell}(T)$ defined by
\[
\mathcal L_T(\Grad r)
\defi
\sum_{F\in\f_T}
\varepsilon_{T,F}(w_{F,n},r)_F
-
(q_h,r)_T,
\qquad
r\in\Poly^{\ell}(T).
\]
This functional is well defined. Indeed, if $r$ is constant, then \eqref{eq:relative.cochain.identity} and \eqref{eq:face.lifting.definition} imply that $\mathcal L_T(\Grad r)=0$.
Let $\vec w_T\in\Grad\Poly^{\ell}(T)\subset\bPoly^{\ell}(T)$ be the $L^2(T)$-Riesz representative of $\mathcal L_T$, characterized by
\begin{equation}
\label{eq:cell.lifting.riesz}
(\vec w_T,\Grad r)_T
=
\sum_{F\in\f_T}
\varepsilon_{T,F}(w_{F,n},r)_F
-
(q_h,r)_T
\qquad
\forall r\in\Poly^{\ell}(T).
\end{equation}
Combining~\eqref{eq:cell.lifting.riesz} with~\eqref{eq:def.DT} gives
\[
\DpT\bigl(\vec w_T,(w_{F,n})_{F\in\f_T}\bigr)
=
q_{h\mid T}.
\]
We can therefore assemble the global hybrid field
\[
\underline{\mathbbmsl w}_{q_h}
\defi
\Bigl(
(\vec w_T)_{T\in\t_h},
(w_{F,n})_{F\in\f_h}
\Bigr).
\]
Since $w_{F,n}=0$ on every boundary face, $\underline{\mathbbmsl w}_{q_h}\in\djdivz$, and the cellwise identity above yields
\begin{equation}
\label{eq:discrete.lifting.right.inverse}
\Dph(\underline{\mathbbmsl w}_{q_h})
=
q_h.
\end{equation}
It remains to prove the stability estimate.
Subtracting the cell average of $r$ in the definition of $\mathcal L_T$ and using the polynomial Poincar\'e and trace inequalities \cite[Lemma~1.32, Lemma~2.15]{Di-Pietro.Droniou:20} gives
\[
|\mathcal L_T(\Grad r)|
\lesssim
\left(
h_T\Norm{0,T}{q_h}
+
\Bigl(
\sum_{F\in\f_T}
h_F\Norm{0,F}{w_{F,n}}^2
\Bigr)^{1/2}
\right)
\Norm{0,T}{\Grad r}.
\]
The Riesz representation~\eqref{eq:cell.lifting.riesz} and the preceding estimate yield
\begin{equation}
\label{eq:cell.lifting.bound}
\Norm{0,T}{\vec w_T}
\lesssim
h_T\Norm{0,T}{q_h}
+
\Bigl(
\sum_{F\in\f_T}
h_F\Norm{0,F}{w_{F,n}}^2
\Bigr)^{1/2}.
\end{equation}
Squaring~\eqref{eq:cell.lifting.bound}, summing over $T\in\t_h$, and using~\eqref{eq:face.lifting.bound}, we obtain
\begin{equation}
\label{eq:cell.lifting.global.bound}
\sum_{T\in\t_h}
\Norm{0,T}{\vec w_T}^2
\lesssim
\Norm{0,\Omega}{q_h}^2.
\end{equation}
\noindent Finally, using the discrete trace inequality~\cite[Lemma~1.32]{Di-Pietro.Droniou:20},~\eqref{eq:face.lifting.bound}, and~\eqref{eq:cell.lifting.global.bound} in \eqref{eq:stab.div} gives
\begin{align*}
s_h(\underline{\mathbbmsl w}_{q_h},\underline{\mathbbmsl w}_{q_h})
&=
\sum_{T\in\t_h}
\sum_{F\in\f_T}
h_F
\Norm{0,F}{\vec w_T\cdot\normal_F-w_{F,n}}^2
\\
&\lesssim
\sum_{T\in\t_h}
\Norm{0,T}{\vec w_T}^2
+
\sum_{F\in\f_h}
h_F\Norm{0,F}{w_{F,n}}^2
\\
&\lesssim
\Norm{0,\Omega}{q_h}^2.
\end{align*}
Combining this estimate with \eqref{eq:discrete.lifting.right.inverse}, \eqref{eq:cell.lifting.global.bound}, and the definition of the discrete graph norm~\eqref{eq:def.discrete.norm}, we conclude that
\[
\tnorm{h}{\underline{\mathbbmsl w}_{q_h}}
\lesssim
\Norm{0,\Omega}{q_h}.
\]
Each step of the construction is linear in $q_h$; hence the lifting can be selected linearly in $q_h$, which concludes the proof.
\end{proof}

\subsection{Discrete inf-sup stability and well-posedness}
\label{sec:discrete.stability}
The standard Brezzi theory requires coercivity of $a_h$ on the kernel of the coupling and a uniform inf-sup condition for $b_h$.
The latter follows directly from the right inverse constructed above.
Define
\[
\ker(b_h)
\defi
\left\{
\jh\in\djdivz:
b_h(\jh,q_h)=0
\quad
\forall q_h\in\Ph
\right\}.
\]

\begin{lemma}[Coercivity on $\ker(b_h)$]
\label{lem:coercive.kernel}
The sesquilinear form $a_h$ is uniformly coercive on $\ker(b_h)$. More precisely,
$$
\Re\,a_h(\jh,\jh)
\geq
c_a\,
\tnorm{h}{\jh}^2, \qquad \forall \jh\in\ker(b_h),
$$
where $c_a\defi\min(\rho_\ast,\alpha)>0$.
\end{lemma}

\begin{proof}
Let $\jh\in\ker(b_h)$.
By~\eqref{eq:Dph.zero.mean}, the reconstructed divergence satisfies $\Dph\jh\in\Ph$. It can therefore be used as a test function in the definition of $\ker(b_h)$, yielding
$$
0
=
b_h(\jh,\Dph\jh)
=
-\Norm{0,\Omega}{\Dph\jh}^2.
$$
Consequently,
\begin{equation}
\Dph\jh=0.\label{eq:kernel.divergence.zero}
\end{equation}
Since $G$ is self-adjoint and positive semidefinite,
\[
g_h(\jh,\jh)
=
(\mu G(\jTh),\jTh)_\Omega
\in
\Real_+.
\]
The contribution $\ii\omega g_h(\jh,\jh)$ is therefore purely imaginary, and hence
\[
\Re\,a_h(\jh,\jh)
=
(\rho\jTh,\jTh)_\Omega
+
\alpha s_h(\jh,\jh)
\geq
\rho_\ast\Norm{0,\Omega}{\jTh}^2
+
\alpha s_h(\jh,\jh).
\]
Using~\eqref{eq:kernel.divergence.zero} in the definition of the discrete graph norm~\eqref{eq:def.discrete.norm} proves Lemma~\ref{lem:coercive.kernel}.
\end{proof}

\begin{lemma}[Uniform discrete inf-sup condition]
\label{lem:infsup.bh}
There exists a constant $\beta>0$, independent of $h$, such that
\begin{equation}
\sup_{\wh\in\djdivz\setminus\{0\}}
\frac{|b_h(\wh,q_h)|}{\tnorm{h}{\wh}}
\geq
\beta\Norm{0,\Omega}{q_h}
\qquad
\forall q_h\in\Ph. \label{eq:discrete.infsup}
\end{equation}
\end{lemma}

\begin{proof}
Let $q_h\in\Ph$ be arbitrary, and let $\underline{\mathbbmsl w}_{q_h}\in\djdivz$ denote the corresponding stable discrete lifting given by Lemma~\ref{lem:stable.discrete.lifting}. By construction, this lifting satisfies $\Dph(\underline{\mathbbmsl w}_{q_h})=q_h$, hence
\[
b_h(\underline{\mathbbmsl w}_{q_h},q_h)
=
-(q_h,q_h)_\Omega
=
-\Norm{0,\Omega}{q_h}^2.
\]
Using the stability of the lifting, we infer
\[
\frac{|b_h(\underline{\mathbbmsl w}_{q_h},q_h)|}
{\tnorm{h}{\underline{\mathbbmsl w}_{q_h}}}
\geq
\frac{1}{C}\Norm{0,\Omega}{q_h}.
\]
Taking the supremum over $\wh\in\djdivz\setminus\{0\}$ concludes the proof with $\beta\defi C^{-1}$.
\end{proof}

\begin{proof}[Proof of Theorem~\ref{thm:well.posedness.discrete}]
The continuity estimates~\eqref{eq:continuity.ah}--\eqref{eq:continuity.bh}, the kernel coercivity of Lemma~\ref{lem:coercive.kernel}, and the uniform inf-sup condition of Lemma~\ref{lem:infsup.bh} satisfy the hypotheses of Brezzi theory in complex Hilbert spaces. Existence and uniqueness of the discrete solution therefore follow.

\noindent It remains to estimate the right-hand side. For every $\wh\in\djdivz$, the Cauchy--Schwarz inequality gives
\[
|\ell_h(\wh)|
=
\omega|(\vec a_0,\wTh)_\Omega|
\leq
\omega\Norm{0,\Omega}{\vec a_0}
\Norm{0,\Omega}{\wTh}
\leq
\omega\Norm{0,\Omega}{\vec a_0}
\tnorm{h}{\wh}.
\]
Hence,
\[
\sup_{\wh\in\djdivz\setminus\{0\}}
\frac{|\ell_h(\wh)|}{\tnorm{h}{\wh}}
\leq
\omega\Norm{0,\Omega}{\vec a_0}.
\]
The stated stability estimate~\eqref{eq:discrete.solution.stability} then follows from the abstract a priori bound associated with Brezzi theory.
\end{proof}

\subsection{Estimate of the consistency error}
\label{sec:consistency}
We insert the interpolated continuous solution into the discrete equations and estimate the resulting residual.
The bound separates the current-projection error, the scalar-potential consistency terms, and the stabilization contribution; the nonlocal operator enters through its $L^2$ boundedness.

Let $\mathcal X_h\defi\djdivz\times\Ph$. We write the mixed formulation in terms of the sesquilinear form
\[
z_h\bigl((\jh,u_h),(\wh,v_h)\bigr)
\defi a_h(\jh,\wh)+\overline{b_h(\wh,u_h)}+b_h(\jh,v_h).
\]
For a conjugate-linear functional $\mathcal R$ on $\mathcal X_h$, its dual norm is
\[
\tnorm{h}{\mathcal R}^{*}
\defi\sup_{(\wh,v_h)\in\mathcal X_h\setminus\{(0,0)\}}
\frac{|\mathcal R(\wh,v_h)|}{\tnorm{h}{(\wh,v_h)}}.
\]

\begin{lemma}[Estimate of the consistency error]\label{lem:consistency}
Let $(\vec j,u)\in\Hzdiv\times\LLstar$ solve Problem~\eqref{eq:weak.model}, and assume that $u\in H^1(\Omega)$ and $\vec j|_T\in\vec H^1(T)$ for every $T\in\t_h$. Define $\mathcal E((\vec j,u);\bullet):\mathcal X_h\to\Comp$ by
\begin{equation}\label{eq:consist.err.def}
\mathcal E\bigl((\vec j,u);(\wh,v_h)\bigr)
\defi\ell_h(\wh)-z_h\bigl((\Ihh\vec j,\proj{h}u),(\wh,v_h)\bigr).
\end{equation}
Under the assumptions of Theorem~\ref{thm:well.posedness.discrete},
\begin{align}\label{eq:consist.err.est}
\tnorm{h}{\mathcal E((\vec j,u);\bullet)}^{*}
&\lesssim\Norm{0,\Omega}{\vec j-\vproj{h}\vec j}\notag\\
&\qquad+\left(\sum_{T\in\t_h}\left[h_T^{-2}\Norm{0,T}{u-\proj{T}u}^{2}
+\Norm{0,T}{\Grad(u-\proj{T}u)}^{2}\right]\right)^{1/2}\notag\\
&\qquad+\left(\sum_{T\in\t_h}\sum_{F\in\f_T}
h_F\Norm{0,F}{(\vec j-\vproj{T}\vec j)\cdot\normal_F}^{2}\right)^{1/2}.
\end{align}
The hidden constant is independent of $h$.
\end{lemma}

\begin{proof}
Let $(\wh,v_h)\in\mathcal X_h$ and introduce the bounded operator
\[
K:\bLL\to\bLL,\qquad K\vec v=\rho\vec v+\ii\omega\mu G(\vec v).
\]
Set $\delta\vec j=\vec j-\vproj{h}\vec j$, $\delta_Tu=u-\proj{T}u$, and $d_{T,F}(\wh)=\wT\cdot\normal_F-\wFn$. The first equation of~\eqref{eq:strong.model}, the definition of $z_h$, and the commuting property~\eqref{eq:global.commuting.Dph} give
\begin{align*}
\mathcal E\bigl((\vec j,u);(\wh,v_h)\bigr)
&=(K\delta\vec j,\wTh)_\Omega+(\Grad u,\wTh)_\Omega
+(\proj{h}u,\Dph\wh)_\Omega\\
&\qquad-\alpha s_h(\Ihh\vec j,\wh),
\end{align*}
where $\Dph(\Ihh\vec j)=\proj{h}(\Div\vec j)=0$ eliminates the contribution involving $v_h$.
Since $u\in H^1(\Omega)$, its traces coincide on interior faces. The single-valued face unknowns and the homogeneous boundary condition therefore imply
\begin{equation}\label{eq:vanishing.jumps}
\sum_{T\in\t_h}\sum_{F\in\f_T}\varepsilon_{T,F}(u,w_{F,n})_F=0.
\end{equation}
Using the definition of $\Dph$, integrating by parts on each cell, and noting that $(\delta_Tu,\Div\wT)_T=0$ by polynomial orthogonality, we obtain
\[
(\Grad u,\wTh)_\Omega+(\proj{h}u,\Dph\wh)_\Omega
=\sum_{T\in\t_h}\sum_{F\in\f_T}\varepsilon_{T,F}(\delta_Tu,d_{T,F}(\wh))_F.
\]
Consequently,
\begin{align}\label{eq:consistency.decomposition}
\mathcal E\bigl((\vec j,u);(\wh,v_h)\bigr)
&=(K\delta\vec j,\wTh)_\Omega
+\sum_{T\in\t_h}\sum_{F\in\f_T}\varepsilon_{T,F}(\delta_Tu,d_{T,F}(\wh))_F\notag\\
&\qquad-\alpha s_h(\Ihh\vec j,\wh).
\end{align}
If the resistivity is constant on each cell, $\rho|_T=\rho_T$, polynomial orthogonality gives
\[
(\rho\delta\vec j,\wTh)_\Omega
=\sum_{T\in\t_h}\rho_T
(\vec j-\vproj{T}\vec j,\wT)_T=0.
\]
For cellwise constant resistivity, the first term in~\eqref{eq:consistency.decomposition} therefore reduces to the Green contribution. The estimate below also covers resistivity varying within cells.

\noindent We estimate the three terms separately. Boundedness of $K$ yields
\[
|(K\delta\vec j,\wTh)_\Omega|
\lesssim\Norm{0,\Omega}{\delta\vec j}\,\tnorm{h}{\wh}.
\]
For the scalar-potential contribution, the weighted Cauchy--Schwarz inequality gives
\begin{align*}
\left|\sum_{T\in\t_h}\sum_{F\in\f_T}\varepsilon_{T,F}(\delta_Tu,d_{T,F}(\wh))_F\right|
&\leq\left(\sum_{T\in\t_h}\sum_{F\in\f_T}h_F^{-1}\Norm{0,F}{\delta_Tu}^{2}\right)^{1/2}
s_h(\wh,\wh)^{1/2}.
\end{align*}
The continuous trace inequality~\cite[Lemma~1.31]{Di-Pietro.Droniou:20} implies
\[
\Norm{0,F}{\delta_Tu}^{2}
\lesssim h_T^{-1}\Norm{0,T}{\delta_Tu}^{2}+h_T\Norm{0,T}{\Grad\delta_Tu}^{2}.
\]
Using $h_F\simeq h_T$, the uniformly bounded number of faces per cell, and~\eqref{eq:def.discrete.norm}, we infer
\[
\left|\sum_{T\in\t_h}\sum_{F\in\f_T}\varepsilon_{T,F}(\delta_Tu,d_{T,F}(\wh))_F\right|
\lesssim\left(\sum_{T\in\t_h}\left[h_T^{-2}\Norm{0,T}{\delta_Tu}^{2}
+\Norm{0,T}{\Grad\delta_Tu}^{2}\right]\right)^{1/2}\tnorm{h}{\wh}.
\]
Finally, since $(\vproj{T}\vec j)\cdot\normal_F\in\Poly^\ell(F)$,
\[
(\vproj{T}\vec j)\cdot\normal_F-\proj{F}(\vec j\cdot\normal_F)
=-\proj{F}\bigl((\vec j-\vproj{T}\vec j)\cdot\normal_F\bigr).
\]
The $L^2(F)$ contraction property of $\proj{F}$ and the Cauchy--Schwarz inequality for $s_h$ therefore give
\begin{align*}
|\alpha s_h(\Ihh\vec j,\wh)|
&\leq\alpha s_h(\Ihh\vec j,\Ihh\vec j)^{1/2}s_h(\wh,\wh)^{1/2}\\
&\lesssim\left(\sum_{T\in\t_h}\sum_{F\in\f_T}h_F
\Norm{0,F}{(\vec j-\vproj{T}\vec j)\cdot\normal_F}^{2}\right)^{1/2}\tnorm{h}{\wh}.
\end{align*}
Combining these estimates in~\eqref{eq:consistency.decomposition}, using $\tnorm{h}{\wh}\leq\tnorm{h}{(\wh,v_h)}$, and taking the supremum over nonzero test pairs proves~\eqref{eq:consist.err.est}.
\end{proof}

\subsection{Proof of the a priori error estimate}
The consistency bound and uniform mixed stability yield Theorem~\ref{thm:error.estimate}.

\begin{proof}[Proof of Theorem~\ref{thm:error.estimate}]
Set $(\underline{\mathbbmsl e}_h,e_h)\defi(\jh,u_h)-(\Ihh\vec j,\proj{h}u)$. The discrete problem and~\eqref{eq:consist.err.def} imply
\[
z_h\bigl((\underline{\mathbbmsl e}_h,e_h),(\wh,v_h)\bigr)
=\mathcal E\bigl((\vec j,u);(\wh,v_h)\bigr)
\qquad\forall(\wh,v_h)\in\mathcal X_h.
\]
The uniform stability of the mixed formulation, established in Theorem~\ref{thm:well.posedness.discrete}, yields
\[
\tnorm{h}{(\underline{\mathbbmsl e}_h,e_h)}
\lesssim\tnorm{h}{\mathcal E((\vec j,u);\bullet)}^{*}.
\]
Lemma~\ref{lem:consistency} proves~\eqref{eq:error.estimate}.
Under the additional regularity assumptions, the approximation properties of the $L^2$ projectors~\cite[Theorem~1.45]{Di-Pietro.Droniou:20} give
\begin{align*}
\Norm{0,T}{\vec j-\vproj{T}\vec j}&\lesssim h_T^{r+1}|\vec j|_{\vec H^{r+1}(T)},\\
h_T^{-1}\Norm{0,T}{u-\proj{T}u}+\Norm{0,T}{\Grad(u-\proj{T}u)}
&\lesssim h_T^r|u|_{H^{r+1}(T)}.
\end{align*}
The continuous trace inequality and the same approximation estimates also imply
\[
\sum_{F\in\f_T}h_F\Norm{0,F}{(\vec j-\vproj{T}\vec j)\cdot\normal_F}^{2}
\lesssim h_T^{2(r+1)}|\vec j|_{\vec H^{r+1}(T)}^2.
\]
Substitution into~\eqref{eq:error.estimate} proves~\eqref{eq:error.rate}.
\end{proof}

\section{Conclusions}
\label{sec:conclusions}
The proposed HHO method combines arbitrary-order approximation with unstructured polyhedral meshes for three-dimensional volume integral eddy-current simulations.
Local refinement allows the mesh to follow a defect boundary while retaining larger cells elsewhere. On any such mesh the polynomial degree can also be increased without changing the cell geometry.
For bounded, uniformly positive resistivity and uniform magnetic permeability, we establish stability and an a priori error estimate for the current.

\noindent For the smooth test problem, selected coarse degree-two discretizations reduced both the cell-current error and the number of unknowns in the dense Green block compared with finer degree-one discretizations. The savings in storage did not always translate into shorter assembly times.

\noindent The defective-plate benchmark illustrates the complementary role of local refinement. At degree one, the relative $\ell^2$ discrepancy between the computed and measured reference-subtracted complex impedance scans was $2.94\%$ on the finest mesh.
The air-cored probe can be moved over the plate without remeshing. Since only the conductor is discretized, the mesh refined around the slot remains fixed and there is no air mesh to deform or regenerate as the coil moves.
At fixed frequency, only the prescribed excitation and the impedance evaluation change with the probe position. We therefore assembled and factorized the system matrix just once for all 76 positions.

\noindent The radial slot provides a benchmark with a known defect geometry. For less regular defects, repeated forward simulations could help assess whether changes in defect shape and local resistivity produce distinguishable impedance signatures. Comparing these predictions with measurements would require separating discretization error from uncertainty in the defect geometry, material properties and experimental data.

\noindent The assembly and factorization costs point to two priorities for larger problems. Compression and fast summation could reduce the cost of distant Green interactions, while preconditioned iterative solvers could reduce the cost of solving the mixed system. For nearby cells, efficient quadrature remains important.
Error estimators for the impedance could guide mesh refinement and the choice of polynomial degree, possibly combined with singularity-adapted enrichment near defect tips.

\noindent A formulation based on a vector potential for the current could enforce the divergence constraint through a compatible discrete curl~\cite{Di-Pietro.Droniou:23DDR,Pitassi:26}. Once the gauge and topological degrees of freedom are accounted for, its system size and computational cost could be compared directly with those of the present mixed formulation.
Anisotropic resistivity could be incorporated through the local resistive term, whereas spatially varying permeability would require the integral representation to account for an additional magnetic response.

\section*{Acknowledgements}
The authors thank Professor J\'er\^ome Droniou for making the HArDCore3D library available for the numerical implementation of this work.
The work of S. Pitassi was supported by the European Union's Horizon Europe research and innovation programme under the Marie Sk\l{}odowska-Curie grant agreement No.~101146324 (EffECT), and by the Japan Society for the Promotion of Science through the JSPS Postdoctoral Fellowship for Research in Japan (Standard), fellowship No.~P26320.
The work of R. Kumar and A. Spadotto was funded by the European Union (ERC Synergy, NEMESIS, project number 101115663).
Views and opinions expressed are those of the authors only and do not necessarily reflect those of the European Union or the European Research Council Executive Agency.
Neither the European Union nor the granting authority can be held responsible for them.

\clearpage
\appendix
\section{Excitation and mesh data for the ECT benchmark}
\label{app:ect.setup}
For a coil centred at $\vec x_{\rm c}=(x_{\rm c},y_{\rm c})$, let
\begin{align*}
r_{\rm c}(x,y)&=\sqrt{(x-x_{\rm c})^2+(y-y_{\rm c})^2},\\
\Omega_{\rm c}(\vec x_{\rm c})&=\{(x,y,z):r_{\rm i}\leq r_{\rm c}\leq r_{\rm o},\ l_0\leq z\leq l_0+h_{\rm c}\}.
\end{align*}
The homogenized winding current and its excitation potential are
\begin{align*}
\vec J_{\rm c}^{\,\vec x_{\rm c}}(x,y,z)
&=\frac{NI}{(r_{\rm o}-r_{\rm i})h_{\rm c}}\,
\frac{(-(y-y_{\rm c}),\,x-x_{\rm c},\,0)^{\mathsf T}}{r_{\rm c}(x,y)},\\
\vec a_0^{\,\vec x_{\rm c}}(\vec x)
&=\frac{\mu_0}{4\pi}\int_{\Omega_{\rm c}(\vec x_{\rm c})}
\frac{\vec J_{\rm c}^{\,\vec x_{\rm c}}(\vec y)}{|\vec x-\vec y|}\,\mathrm d\vec y.
\end{align*}
The winding-current formula applies inside $\Omega_{\rm c}(\vec x_{\rm c})$ and is extended by zero outside the coil. The positive lift-off makes this source integral regular on the conductor. The coil quadrature used for the scans is specified in Section~\ref{sec:numerics.fastener.hole}. Cell integration of the right-hand side and of~\eqref{eq:reaction.impedance} uses polynomial quadrature degrees $2\ell+6$ and $2\ell+8$, respectively.

\noindent For the Green operator, distinct-cell interactions use product quadrature with parameter $q_G$.
For self-interactions, the inner integration uses an oriented radial formula.
For a fixed outer quadrature point $\vec x$ and the polynomial extension of $\vec p\in\bPoly^\ell(T)$ to $\Real^3$, this formula reads
\begin{equation}
\label{eq:green.radial.identity}
\begin{aligned}
\int_T\frac{\vec p(\vec y)}{4\pi|\vec y-\vec x|}\,\mathrm d\vec y
&=\frac{1}{4\pi}\sum_{F\in\f_T}\int_F
\frac{(\vec z-\vec x)\cdot\normal_{T,F}}{|\vec z-\vec x|}
\\[-0.3em]
&\qquad\times\left[\int_0^1 t\,\vec p\bigl(\vec x+t(\vec z-\vec x)\bigr)\,\mathrm dt\right]\,\mathrm dS_{\vec z}.
\end{aligned}
\end{equation}
The identity follows componentwise from the divergence theorem and retains the signed face factors; it does not require a partition of the cell into pyramids contained in $T$.
The transformation $\vec y=\vec x+t(\vec z-\vec x)$ combines the $1/t$ kernel singularity with the $t^2$ volume Jacobian, leaving the factor $t$ in~\eqref{eq:green.radial.identity}.
The self-interaction quadrature parameter is increased to $q_G+4$.
The sensitivity checks in Section~\ref{sec:numerics.fastener.hole} compare $q_G=4$ and $5$.

\noindent The meshes are supplied in the HArDCore3D RF format as paired \texttt{.node} and \texttt{.ele} files, with coordinates in metres. The hole boundary uses $128$ chords, with a maximum radial deviation of approximately $0.003$~mm, held fixed across levels. The slot has planar sides at $y=\pm0.117$~mm and a flat end at $x=19.8$~mm. Table~\ref{tab:ect.mesh.details} reports the face counts and cell diameters; $h_{\min}$ is a cell diameter and is distinct from the slot width.

\begin{table}[tb]
\centering
\caption{Additional ECT mesh data. Cell diameters and layer thicknesses are in millimetres.}
\label{tab:ect.mesh.details}
\small
\begin{tabular}{lrrrrrr}
\hline
Mesh & $N_F$ & Interior faces & Layers & Layer thickness & $h_{\min}$ & $h_{\max}$ \\
\hline
L0 & 1612 & 844 & 1 & 2.000 & 2.366 & 108.317 \\
L1 & 4585 & 3235 & 2 & 1.000 & 1.311 & 82.804 \\
L2 & 13063 & 10478 & 3 & 0.667 & 0.803 & 55.192 \\
\hline
\end{tabular}
\end{table}

\FloatBarrier
\section{Proof of the relative cochain lifting}
\label{app:cochain.proof}
We give the proof of Lemma~\ref{lem:relative.cochain.poincare}, retaining the relative boundary condition required by the homogeneous normal-current constraint.

\begin{proof}
The argument is a relative, top-degree version of the construction used in the proof of~\cite[Lemma~6]{DiPietro.Droniou.Hanot.Pitassi:25}, expressed through the Raviart--Thomas vector proxy.
Let $\mathcal S_h$ be the conforming simplicial submesh associated with $\t_h$. Let $\widehat{\f}_h$ denote the set of two-dimensional simplices of $\mathcal S_h$, and, for every $F\in\f_h$, set
\[
\widehat{\f}_h(F)
\defi
\left\{
\widehat F\in\widehat{\f}_h
\mid
\widehat F\subset F
\right\}.
\]
We define the piecewise constant function $\widetilde\xi_h$ on $\mathcal S_h$ by
\begin{equation}
\label{eq:def.relative.simplicial.data}
\widetilde\xi_{h\mid S}
\defi
\frac{\xi_T}{|T|}
\qquad
\forall S\in\mathcal S_h,
\quad
S\subset T.
\end{equation}
By the compatibility condition~\eqref{eq:relative.cochain.compatibility}, we have
$$
\int_\Omega\widetilde\xi_h
=
\sum_{T\in\t_h}
\sum_{S\subset T}
|S|\frac{\xi_T}{|T|}
=
\sum_{T\in\t_h}\xi_T
=
0.
$$
Now, consider the lowest-order Raviart--Thomas space with homogeneous normal trace
$$
\boldsymbol{\mathrm{RT}}_{0,0}(\mathcal S_h)
\defi
\left\{
\vec v_h\in\Hzdiv
\mid
\vec v_{h\mid S}\in\boldsymbol{\mathrm{RT}}_0(S)
\quad
\forall S\in\mathcal S_h
\right\}.
$$
The finite element Poincar\'e inequality for the de Rham complex with homogeneous boundary conditions~\cite[Section~6.2]{Arnold.Falk.Winther:10}, applied to the lowest-order Raviart--Thomas vector proxy, yields $\widetilde{\vec\lambda}_h\in\boldsymbol{\mathrm{RT}}_{0,0}(\mathcal S_h)$ such that
\begin{equation}
\label{eq:relative.RT.lifting}
\Div\widetilde{\vec\lambda}_h
=
\widetilde\xi_h,
\qquad
\Norm{0,\Omega}{\widetilde{\vec\lambda}_h}
\lesssim
\Norm{0,\Omega}{\widetilde\xi_h}.
\end{equation}
This is the homogeneous-normal counterpart of the Raviart--Thomas lifting underlying the divergence Poincar\'e argument in~\cite{DiPietro.Hanot:24}.
It is also the relative counterpart of the finite element Poincar\'e estimate used in the proof of~\cite[Lemma~6]{DiPietro.Droniou.Hanot.Pitassi:25}.
Selecting the solution of minimum $L^2$ norm makes the construction linear in $\widetilde\xi_h$.

\noindent Fix a unit normal $\normal_{\widehat F}$ on every $\widehat F\in\widehat{\f}_h$ and define the corresponding oriented flux by
\[
\widetilde\lambda_{\widehat F}
\defi
\int_{\widehat F}
\widetilde{\vec\lambda}_h\cdot\normal_{\widehat F}.
\]
The fact $\widetilde{\vec\lambda}_h\in\Hzdiv$ leads to
\begin{equation}
\label{eq:relative.simplicial.face.zero}
\widetilde\lambda_{\widehat F}=0
\qquad
\forall\widehat F\subset\partial\Omega.
\end{equation}
For every polyhedral face $F\in\f_h$, we collect the corresponding simplicial fluxes by setting
\begin{equation}
\label{eq:relative.cochain.aggregation}
\lambda_F
\defi
\sum_{\widehat F\in\widehat{\f}_h(F)}
\varepsilon_{F,\widehat F}
\widetilde\lambda_{\widehat F},
\end{equation}
where $\varepsilon_{F,\widehat F}\in\{-1,1\}$ compares the orientation induced by $\normal_{\widehat F}$ with that induced by $\normal_F$.
The boundary condition~\eqref{eq:relative.cochain.boundary} follows from~\eqref{eq:relative.simplicial.face.zero}.

\noindent Let $T\in\t_h$. Summing the divergence theorem over all tetrahedra $S\subset T$, the flux contributions across the internal faces of the simplicial submesh cancel pairwise. Consequently,
\begin{align*}
\sum_{F\in\f_T}
\varepsilon_{T,F}\lambda_F
&=
\sum_{S\subset T}
\int_{\partial S}
\widetilde{\vec\lambda}_h\cdot\normal_S
\\
&=
\sum_{S\subset T}
\int_S
\Div\widetilde{\vec\lambda}_h
\\
&=
\sum_{S\subset T}
|S|\frac{\xi_T}{|T|}
=
\xi_T,
\end{align*}
which proves~\eqref{eq:relative.cochain.identity}.

\noindent It remains to establish the estimate~\eqref{eq:relative.cochain.bound}.
The uniform bound on the cardinality of $\widehat{\f}_h(F)$ for $F\in\f_h$, the Cauchy--Schwarz inequality, local diameter comparability, and~\eqref{eq:relative.cochain.aggregation} give
\[
\sum_{F\in\f_h}
h_F^{-1}|\lambda_F|^2
\lesssim
\sum_{\widehat F\in\widehat{\f}_h}
h_{\widehat F}^{-1}
|\widetilde\lambda_{\widehat F}|^2 \lesssim
 \Norm{0,\Omega}{\widetilde{\vec\lambda}_h}^2,
\]
where the standard scaling of the lowest-order Raviart--Thomas face degrees of freedom is used in the last step.
Combining this estimate with~\eqref{eq:relative.RT.lifting}, we infer that
\[
\sum_{F\in\f_h}
h_F^{-1}|\lambda_F|^2
\lesssim
\Norm{0,\Omega}{\widetilde\xi_h}^2.
\]
Finally, the definition~\eqref{eq:def.relative.simplicial.data} and mesh regularity yield
\begin{align*}
\Norm{0,\Omega}{\widetilde\xi_h}^2
&=
\sum_{T\in\t_h}
\sum_{S\subset T}
|S|\frac{|\xi_T|^2}{|T|^2}
\\
&=
\sum_{T\in\t_h}
\frac{|\xi_T|^2}{|T|}
\lesssim
\sum_{T\in\t_h}
h_T^{-3}|\xi_T|^2.
\end{align*}
This proves~\eqref{eq:relative.cochain.bound} and concludes the proof.
\end{proof}

\bibliographystyle{plainurl}
\bibliography{references}

\end{document}

%% file: figures/ect-current.tex
\begin{figure}[tb]
\centering
\begin{tikzpicture}
\begin{axis}[width=.80\textwidth,height=.40\textwidth,axis equal image,
xmin=7,xmax=26,ymin=-8,ymax=8,xlabel={$x$ (mm)},ylabel={$y$ (mm)},
tick label style={font=\scriptsize},label style={font=\small},colormap/viridis,
colorbar,colorbar style={ylabel={$|\Re\jTh|/(10^7\,\mathrm{A/m^2})$},ylabel style={font=\small},tick label style={font=\scriptsize}}]
\addplot[scatter,only marks,scatter src=explicit,mark=*,mark size=1.25pt] coordinates {
(8.79,5.57648) [3.49172]
(20.3325,-0.360618) [6.53581]
(10.3768,-2.27547) [4.14773]
(10.3081,1.61909) [3.30771]
(9.92825,2.73423) [4.79113]
(19.9474,0.302589) [10.1122]
(21.5182,-1.15638) [0.948735]
(8.37334,-6.56558) [2.13009]
(23.0005,0.725491) [2.4865]
(14.8521,0.50635) [5.37782]
(19.3429,-0.756282) [6.73352]
(9.64088,3.83376) [4.99709]
(7.40561,7.56058) [0.97042]
(10.7399,3.38127) [4.48789]
(9.96455,-4.13362) [4.80619]
(9.87701,-3.15495) [5.04652]
(7.45047,-7.22045) [0.745873]
(6.56552,-7.93348) [1.28648]
(9.29255,-4.69145) [4.57491]
(8.98617,-5.71514) [3.44684]
(11.0358,-0.361628) [2.05252]
(10.3608,-0.89505) [1.97246]
(16.8257,-0.51407) [5.686]
(8.20343,6.4082) [2.06768]
(9.98377,4.81876) [4.27682]
(9.49517,6.57031) [3.10279]
(10.557,0.523191) [1.4751]
(11.7544,0.472238) [3.29273]
(20.213,0.175601) [7.22615]
(10.1548,-5.3555) [3.93817]
(12.1752,-2.24256) [3.75378]
(19.9772,-0.609395) [7.08695]
(23.2656,-4.10343) [3.9662]
(11.5204,-3.80656) [3.68465]
(11.127,2.06493) [3.93445]
(8.55341,7.82471) [2.4045]
(11.0727,-1.84762) [3.63322]
(17.2427,0.640888) [5.83145]
(18.487,-0.955922) [6.00679]
(11.9693,-7.3197) [3.38451]
(17.9712,1.29536) [5.55232]
(12.4,-0.657223) [4.1182]
(12.2564,4.24746) [2.82684]
(20.1414,1.33764) [3.94181]
(11.151,5.85761) [3.32411]
(11.2331,1.15953) [3.17971]
(10.8916,7.99692) [3.5397]
(23.0642,-2.65577) [3.33868]
(13.412,-0.493344) [4.77903]
(22.0736,0.227111) [0.756013]
(21.3681,0.00120491) [1.13886]
(11.3687,-5.26003) [3.34508]
(13.4315,-1.72717) [4.48941]
(25.3005,-2.30152) [4.04155]
(19.4344,0.874063) [6.2476]
(19.9604,-4.05435) [1.79231]
(18.0942,7.50338) [3.84771]
(18.2862,0.360871) [5.86294]
(19.9204,-0.0617988) [11.0618]
(21.4139,-3.67116) [2.63295]
(13.4345,-5.52349) [1.96711]
(19.2128,-1.47002) [4.94755]
(20.3106,-1.33229) [3.59281]
(17.1137,-3.14397) [2.57671]
(20.9009,-0.811755) [2.38673]
(20.3421,0.78068) [4.84071]
(18.012,4.39866) [1.00079]
(18.8601,2.14515) [3.82449]
(13.7412,3.41571) [2.43841]
(17.6277,2.61105) [3.48413]
(20.0896,4.09198) [1.85673]
(20.7514,-0.180735) [3.69238]
(18.8408,0.954789) [6.0758]
(20.9513,0.707547) [2.50621]
(19.3408,0.377616) [7.45244]
(20.7628,-1.94708) [1.86658]
(21.68,-2.19395) [1.75768]
(21.7865,2.01369) [1.60338]
(18.9956,3.22963) [1.77889]
(19.0123,-0.332724) [6.28734]
(18.025,-0.355856) [5.7057]
(18.0844,-1.77027) [4.90034]
(20.3648,-2.61518) [1.69183]
(15.5891,-0.533814) [5.54556]
(17.7166,-1.02005) [5.8978]
(17.1697,-1.59447) [5.25929]
(18.3434,-2.91954) [2.62772]
(11.652,-0.847534) [3.35914]
(19.3582,1.59343) [4.62793]
(16.8104,1.49929) [5.35964]
(20.473,2.10958) [2.03759]
(15.754,1.61239) [5.15849]
(19.2966,-2.89922) [2.14359]
(21.4772,4.33792) [3.17697]
(16.128,0.495634) [5.59873]
(14.3601,5.45663) [1.45722]
(16.0159,5.07685) [0.60519]
(16.4994,3.23045) [2.49655]
(23.8694,-0.744682) [3.37432]
(22.59,-1.53979) [2.25729]
(19.8822,-1.97432) [3.07474]
(15.49,-3.25231) [2.47391]
(13.5774,-3.41414) [2.50051]
(17.7342,-4.6854) [0.964666]
(15.3884,-4.9528) [0.484544]
(12.7862,-1.25388) [4.50833]
(14.9619,-1.48601) [5.20711]
(14.1752,-1.20855) [5.22622]
(15.2023,2.60397) [3.59996]
(12.834,0.878953) [4.493]
(12.2396,2.58342) [3.62434]
(23.141,4.8833) [4.14474]
(25.9855,0.0422744) [3.94612]
(16.2776,-1.21811) [5.73877]
(16.0843,-1.95573) [4.77493]
(24.5602,0.872674) [3.83145]
(24.9208,-5.22398) [3.86687]
(19.0646,-5.88954) [2.98239]
(14.0281,1.1876) [5.23186]
(13.7514,2.03066) [4.32752]
(12.1003,1.35698) [3.96537]
(13.6113,0.432017) [4.88863]
(23.6496,1.90829) [3.52052]
(16.6871,-6.56526) [2.81062]
(13.3347,6.92666) [3.01402]
(22.994,2.9763) [3.44609]
(19.3932,5.88703) [3.21141]
(21.2843,-5.48426) [3.81148]
(25.1826,5.18272) [3.68739]
(22.9636,7.47965) [3.83409]
(23.0127,-6.53344) [4.25857]
(14.6456,-7.37344) [3.3557]
(14.633,-0.445573) [5.30171]
};
\draw[->,black!75,line width=.25pt] (axis cs:8.79,5.57648)--(axis cs:8.70410007,5.81381711);
\draw[->,black!75,line width=.25pt] (axis cs:20.3325,-0.360618)--(axis cs:20.5388035,0.0646566675);
\draw[->,black!75,line width=.25pt] (axis cs:10.3768,-2.27547)--(axis cs:10.4883189,-1.99703073);
\draw[->,black!75,line width=.25pt] (axis cs:10.3081,1.61909)--(axis cs:10.2151459,1.83947218);
\draw[->,black!75,line width=.25pt] (axis cs:9.92825,2.73423)--(axis cs:9.82291381,3.06429414);
\draw[->,black!75,line width=.25pt] (axis cs:19.9474,0.302589)--(axis cs:19.4875241,0.871223317);
\draw[->,black!75,line width=.25pt] (axis cs:21.5182,-1.15638)--(axis cs:21.5844562,-1.17420965);
\draw[->,black!75,line width=.25pt] (axis cs:8.37334,-6.56558)--(axis cs:8.38728637,-6.41234131);
\draw[->,black!75,line width=.25pt] (axis cs:23.0005,0.725491)--(axis cs:22.9951123,0.545745234);
\draw[->,black!75,line width=.25pt] (axis cs:14.8521,0.50635)--(axis cs:14.4632216,0.51255299);
\draw[->,black!75,line width=.25pt] (axis cs:19.3429,-0.756282)--(axis cs:19.8183343,-0.650942197);
\draw[->,black!75,line width=.25pt] (axis cs:9.64088,3.83376)--(axis cs:9.51824247,4.17365906);
\draw[->,black!75,line width=.25pt] (axis cs:7.40561,7.56058)--(axis cs:7.45949815,7.60502513);
\draw[->,black!75,line width=.25pt] (axis cs:10.7399,3.38127)--(axis cs:10.6294377,3.68645211);
\draw[->,black!75,line width=.25pt] (axis cs:9.96455,-4.13362)--(axis cs:10.0652316,-3.80096967);
\draw[->,black!75,line width=.25pt] (axis cs:9.87701,-3.15495)--(axis cs:9.99453661,-2.80945691);
\draw[->,black!75,line width=.25pt] (axis cs:7.45047,-7.22045)--(axis cs:7.42755424,-7.17224599);
\draw[->,black!75,line width=.25pt] (axis cs:6.56552,-7.93348)--(axis cs:6.47794584,-7.96414452);
\draw[->,black!75,line width=.25pt] (axis cs:9.29255,-4.69145)--(axis cs:9.40857016,-4.38167912);
\draw[->,black!75,line width=.25pt] (axis cs:8.98617,-5.71514)--(axis cs:9.0489678,-5.474014);
\draw[->,black!75,line width=.25pt] (axis cs:11.0358,-0.361628)--(axis cs:11.1791061,-0.322945477);
\draw[->,black!75,line width=.25pt] (axis cs:10.3608,-0.89505)--(axis cs:10.4434534,-0.778812876);
\draw[->,black!75,line width=.25pt] (axis cs:16.8257,-0.51407)--(axis cs:17.2369119,-0.513883901);
\draw[->,black!75,line width=.25pt] (axis cs:8.20343,6.4082)--(axis cs:8.1636278,6.55211363);
\draw[->,black!75,line width=.25pt] (axis cs:9.98377,4.81876)--(axis cs:9.91029191,5.1191711);
\draw[->,black!75,line width=.25pt] (axis cs:9.49517,6.57031)--(axis cs:9.53670103,6.79080211);
\draw[->,black!75,line width=.25pt] (axis cs:10.557,0.523191)--(axis cs:10.4684271,0.582648223);
\draw[->,black!75,line width=.25pt] (axis cs:11.7544,0.472238)--(axis cs:11.5190915,0.508808245);
\draw[->,black!75,line width=.25pt] (axis cs:20.213,0.175601)--(axis cs:20.0764566,0.680049027);
\draw[->,black!75,line width=.25pt] (axis cs:10.1548,-5.3555)--(axis cs:10.1800389,-5.0718312);
\draw[->,black!75,line width=.25pt] (axis cs:12.1752,-2.24256)--(axis cs:12.3968883,-2.08586509);
\draw[->,black!75,line width=.25pt] (axis cs:19.9772,-0.609395)--(axis cs:20.4218617,-0.354521508);
\draw[->,black!75,line width=.25pt] (axis cs:23.2656,-4.10343)--(axis cs:23.1515253,-4.36661028);
\draw[->,black!75,line width=.25pt] (axis cs:11.5204,-3.80656)--(axis cs:11.610198,-3.55567139);
\draw[->,black!75,line width=.25pt] (axis cs:11.127,2.06493)--(axis cs:10.9530457,2.29009622);
\draw[->,black!75,line width=.25pt] (axis cs:8.55341,7.82471)--(axis cs:8.66048623,7.96168877);
\draw[->,black!75,line width=.25pt] (axis cs:11.0727,-1.84762)--(axis cs:11.235791,-1.64162028);
\draw[->,black!75,line width=.25pt] (axis cs:17.2427,0.640888)--(axis cs:16.8209734,0.642789188);
\draw[->,black!75,line width=.25pt] (axis cs:18.487,-0.955922)--(axis cs:18.9208736,-0.934349745);
\draw[->,black!75,line width=.25pt] (axis cs:11.9693,-7.3197)--(axis cs:11.7835726,-7.16026919);
\draw[->,black!75,line width=.25pt] (axis cs:17.9712,1.29536)--(axis cs:17.5696603,1.2976452);
\draw[->,black!75,line width=.25pt] (axis cs:12.4,-0.657223)--(axis cs:12.6949834,-0.616125463);
\draw[->,black!75,line width=.25pt] (axis cs:12.2564,4.24746)--(axis cs:12.2019029,4.44450132);
\draw[->,black!75,line width=.25pt] (axis cs:20.1414,1.33764)--(axis cs:19.8705024,1.42642192);
\draw[->,black!75,line width=.25pt] (axis cs:11.151,5.85761)--(axis cs:11.2143761,6.08950713);
\draw[->,black!75,line width=.25pt] (axis cs:11.2331,1.15953)--(axis cs:11.0431678,1.28915807);
\draw[->,black!75,line width=.25pt] (axis cs:10.8916,7.99692)--(axis cs:11.0885126,8.16049986);
\draw[->,black!75,line width=.25pt] (axis cs:23.0642,-2.65577)--(axis cs:23.0116361,-2.89143577);
\draw[->,black!75,line width=.25pt] (axis cs:13.412,-0.493344)--(axis cs:13.7572725,-0.477771419);
\draw[->,black!75,line width=.25pt] (axis cs:22.0736,0.227111)--(axis cs:22.0667975,0.17286038);
\draw[->,black!75,line width=.25pt] (axis cs:21.3681,0.00120491)--(axis cs:21.3591486,0.0830801806);
\draw[->,black!75,line width=.25pt] (axis cs:11.3687,-5.26003)--(axis cs:11.3419055,-5.01960117);
\draw[->,black!75,line width=.25pt] (axis cs:13.4315,-1.72717)--(axis cs:13.7485652,-1.65727383);
\draw[->,black!75,line width=.25pt] (axis cs:25.3005,-2.30152)--(axis cs:25.2419552,-2.58788493);
\draw[->,black!75,line width=.25pt] (axis cs:19.4344,0.874063)--(axis cs:18.9948492,0.978654281);
\draw[->,black!75,line width=.25pt] (axis cs:19.9604,-4.05435)--(axis cs:19.9397067,-4.18230963);
\draw[->,black!75,line width=.25pt] (axis cs:18.0942,7.50338)--(axis cs:18.3504671,7.39492895);
\draw[->,black!75,line width=.25pt] (axis cs:18.2862,0.360871)--(axis cs:17.8622683,0.367839058);
\draw[->,black!75,line width=.25pt] (axis cs:19.9204,-0.0617988)--(axis cs:19.9524877,0.737557427);
\draw[->,black!75,line width=.25pt] (axis cs:21.4139,-3.67116)--(axis cs:21.3613223,-3.85417467);
\draw[->,black!75,line width=.25pt] (axis cs:13.4345,-5.52349)--(axis cs:13.3570985,-5.40412572);
\draw[->,black!75,line width=.25pt] (axis cs:19.2128,-1.47002)--(axis cs:19.5673267,-1.4216618);
\draw[->,black!75,line width=.25pt] (axis cs:20.3106,-1.33229)--(axis cs:20.5548929,-1.24377567);
\draw[->,black!75,line width=.25pt] (axis cs:17.1137,-3.14397)--(axis cs:17.2989393,-3.16428445);
\draw[->,black!75,line width=.25pt] (axis cs:20.9009,-0.811755)--(axis cs:21.0231969,-0.689943378);
\draw[->,black!75,line width=.25pt] (axis cs:20.3421,0.78068)--(axis cs:20.0672732,0.997540877);
\draw[->,black!75,line width=.25pt] (axis cs:18.012,4.39866)--(axis cs:18.0052972,4.3265929);
\draw[->,black!75,line width=.25pt] (axis cs:18.8601,2.14515)--(axis cs:18.5835375,2.1412582);
\draw[->,black!75,line width=.25pt] (axis cs:13.7412,3.41571)--(axis cs:13.5942301,3.51317044);
\draw[->,black!75,line width=.25pt] (axis cs:17.6277,2.61105)--(axis cs:17.3762364,2.59499677);
\draw[->,black!75,line width=.25pt] (axis cs:20.0896,4.09198)--(axis cs:20.1170683,3.96053957);
\draw[->,black!75,line width=.25pt] (axis cs:20.7514,-0.180735)--(axis cs:20.8058931,0.080681368);
\draw[->,black!75,line width=.25pt] (axis cs:18.8408,0.954789)--(axis cs:18.4029277,0.99139678);
\draw[->,black!75,line width=.25pt] (axis cs:20.9513,0.707547)--(axis cs:20.8465822,0.855486735);
\draw[->,black!75,line width=.25pt] (axis cs:19.3408,0.377616)--(axis cs:18.8048721,0.434429747);
\draw[->,black!75,line width=.25pt] (axis cs:20.7628,-1.94708)--(axis cs:20.8968722,-1.96282138);
\draw[->,black!75,line width=.25pt] (axis cs:21.68,-2.19395)--(axis cs:21.7074593,-2.31806577);
\draw[->,black!75,line width=.25pt] (axis cs:21.7865,2.01369)--(axis cs:21.7600172,1.90079627);
\draw[->,black!75,line width=.25pt] (axis cs:18.9956,3.22963)--(axis cs:18.8808787,3.17140578);
\draw[->,black!75,line width=.25pt] (axis cs:19.0123,-0.332724)--(axis cs:19.4668133,-0.321116126);
\draw[->,black!75,line width=.25pt] (axis cs:18.025,-0.355856)--(axis cs:18.4376192,-0.35470554);
\draw[->,black!75,line width=.25pt] (axis cs:18.0844,-1.77027)--(axis cs:18.4387951,-1.77112828);
\draw[->,black!75,line width=.25pt] (axis cs:20.3648,-2.61518)--(axis cs:20.4702902,-2.67716679);
\draw[->,black!75,line width=.25pt] (axis cs:15.5891,-0.533814)--(axis cs:15.9901399,-0.52995755);
\draw[->,black!75,line width=.25pt] (axis cs:17.7166,-1.02005)--(axis cs:18.143128,-1.01851027);
\draw[->,black!75,line width=.25pt] (axis cs:17.1697,-1.59447)--(axis cs:17.5500547,-1.59555);
\draw[->,black!75,line width=.25pt] (axis cs:18.3434,-2.91954)--(axis cs:18.5301298,-2.95485114);
\draw[->,black!75,line width=.25pt] (axis cs:11.652,-0.847534)--(axis cs:11.8819011,-0.769038287);
\draw[->,black!75,line width=.25pt] (axis cs:19.3582,1.59343)--(axis cs:19.0264732,1.63790543);
\draw[->,black!75,line width=.25pt] (axis cs:16.8104,1.49929)--(axis cs:16.4227886,1.49810872);
\draw[->,black!75,line width=.25pt] (axis cs:20.473,2.10958)--(axis cs:20.3278584,2.08410957);
\draw[->,black!75,line width=.25pt] (axis cs:15.754,1.61239)--(axis cs:15.381111,1.62390509);
\draw[->,black!75,line width=.25pt] (axis cs:19.2966,-2.89922)--(axis cs:19.4418472,-2.95340785);
\draw[->,black!75,line width=.25pt] (axis cs:21.4772,4.33792)--(axis cs:21.5756252,4.13030825);
\draw[->,black!75,line width=.25pt] (axis cs:16.128,0.495634)--(axis cs:15.7231039,0.497873368);
\draw[->,black!75,line width=.25pt] (axis cs:14.3601,5.45663)--(axis cs:14.4370423,5.5286472);
\draw[->,black!75,line width=.25pt] (axis cs:16.0159,5.07685)--(axis cs:16.0592652,5.07092679);
\draw[->,black!75,line width=.25pt] (axis cs:16.4994,3.23045)--(axis cs:16.3188479,3.23016449);
\draw[->,black!75,line width=.25pt] (axis cs:23.8694,-0.744682)--(axis cs:23.8570623,-0.988404324);
\draw[->,black!75,line width=.25pt] (axis cs:22.59,-1.53979)--(axis cs:22.5809863,-1.70278985);
\draw[->,black!75,line width=.25pt] (axis cs:19.8822,-1.97432)--(axis cs:20.1042636,-1.96270287);
\draw[->,black!75,line width=.25pt] (axis cs:15.49,-3.25231)--(axis cs:15.667268,-3.22808595);
\draw[->,black!75,line width=.25pt] (axis cs:13.5774,-3.41414)--(axis cs:13.7232048,-3.3071643);
\draw[->,black!75,line width=.25pt] (axis cs:17.7342,-4.6854)--(axis cs:17.7086063,-4.7503014);
\draw[->,black!75,line width=.25pt] (axis cs:15.3884,-4.9528)--(axis cs:15.3620913,-4.9296518);
\draw[->,black!75,line width=.25pt] (axis cs:12.7862,-1.25388)--(axis cs:13.1048468,-1.18482101);
\draw[->,black!75,line width=.25pt] (axis cs:14.9619,-1.48601)--(axis cs:15.3378149,-1.46358364);
\draw[->,black!75,line width=.25pt] (axis cs:14.1752,-1.20855)--(axis cs:14.5517434,-1.17579085);
\draw[->,black!75,line width=.25pt] (axis cs:15.2023,2.60397)--(axis cs:14.9434105,2.63152529);
\draw[->,black!75,line width=.25pt] (axis cs:12.834,0.878953)--(axis cs:12.5124386,0.92567067);
\draw[->,black!75,line width=.25pt] (axis cs:12.2396,2.58342)--(axis cs:12.040977,2.75445471);
\draw[->,black!75,line width=.25pt] (axis cs:23.141,4.8833)--(axis cs:23.2838281,4.61976463);
\draw[->,black!75,line width=.25pt] (axis cs:25.9855,0.0422744)--(axis cs:25.9838919,-0.243107684);
\draw[->,black!75,line width=.25pt] (axis cs:16.2776,-1.21811)--(axis cs:16.6926065,-1.21339922);
\draw[->,black!75,line width=.25pt] (axis cs:16.0843,-1.95573)--(axis cs:16.4295392,-1.9479461);
\draw[->,black!75,line width=.25pt] (axis cs:24.5602,0.872674)--(axis cs:24.5817955,0.596422436);
\draw[->,black!75,line width=.25pt] (axis cs:24.9208,-5.22398)--(axis cs:24.7933423,-5.47290148);
\draw[->,black!75,line width=.25pt] (axis cs:19.0646,-5.88954)--(axis cs:18.8994211,-6.02824001);
\draw[->,black!75,line width=.25pt] (axis cs:14.0281,1.1876)--(axis cs:13.6512594,1.22161956);
\draw[->,black!75,line width=.25pt] (axis cs:13.7514,2.03066)--(axis cs:13.4466113,2.10175345);
\draw[->,black!75,line width=.25pt] (axis cs:12.1003,1.35698)--(axis cs:11.8330538,1.46100212);
\draw[->,black!75,line width=.25pt] (axis cs:13.6113,0.432017)--(axis cs:13.257895,0.442138462);
\draw[->,black!75,line width=.25pt] (axis cs:23.6496,1.90829)--(axis cs:23.6913581,1.65713019);
\draw[->,black!75,line width=.25pt] (axis cs:16.6871,-6.56526)--(axis cs:16.489525,-6.61302459);
\draw[->,black!75,line width=.25pt] (axis cs:13.3347,6.92666)--(axis cs:13.5260612,7.03103721);
\draw[->,black!75,line width=.25pt] (axis cs:22.994,2.9763)--(axis cs:23.0555286,2.73478997);
\draw[->,black!75,line width=.25pt] (axis cs:19.3932,5.88703)--(axis cs:19.5697868,5.73617139);
\draw[->,black!75,line width=.25pt] (axis cs:21.2843,-5.48426)--(axis cs:21.1162702,-5.70277486);
\draw[->,black!75,line width=.25pt] (axis cs:25.1826,5.18272)--(axis cs:25.2965771,4.94162871);
\draw[->,black!75,line width=.25pt] (axis cs:22.9636,7.47965)--(axis cs:23.1494005,7.27382168);
\draw[->,black!75,line width=.25pt] (axis cs:23.0127,-6.53344)--(axis cs:22.8189688,-6.77285995);
\draw[->,black!75,line width=.25pt] (axis cs:14.6456,-7.37344)--(axis cs:14.4057411,-7.33649571);
\draw[->,black!75,line width=.25pt] (axis cs:14.633,-0.445573)--(axis cs:15.0163929,-0.440666494);
\path[fill=white,draw=black,line width=.5pt] (axis cs:-10,-1.2246468e-15)--(axis cs:-9.98795456,-0.490676743)--(axis cs:-9.95184727,-0.980171403)--(axis cs:-9.8917651,-1.46730474)--(axis cs:-9.8078528,-1.95090322)--(axis cs:-9.70031253,-2.4298018)--(axis cs:-9.56940336,-2.90284677)--(axis cs:-9.41544065,-3.36889853)--(axis cs:-9.23879533,-3.82683432)--(axis cs:-9.03989293,-4.27555093)--(axis cs:-8.81921264,-4.71396737)--(axis cs:-8.5772861,-5.14102744)--(axis cs:-8.31469612,-5.55570233)--(axis cs:-8.03207531,-5.95699304)--(axis cs:-7.73010453,-6.34393284)--(axis cs:-7.40951125,-6.71558955)--(axis cs:-7.07106781,-7.07106781)--(axis cs:-6.71558955,-7.40951125)--(axis cs:-6.34393284,-7.73010453)--(axis cs:-5.95699304,-8.03207531)--(axis cs:-5.55570233,-8.31469612)--(axis cs:-5.14102744,-8.5772861)--(axis cs:-4.71396737,-8.81921264)--(axis cs:-4.27555093,-9.03989293)--(axis cs:-3.82683432,-9.23879533)--(axis cs:-3.36889853,-9.41544065)--(axis cs:-2.90284677,-9.56940336)--(axis cs:-2.4298018,-9.70031253)--(axis cs:-1.95090322,-9.8078528)--(axis cs:-1.46730474,-9.8917651)--(axis cs:-0.980171403,-9.95184727)--(axis cs:-0.490676743,-9.98795456)--(axis cs:6.123234e-16,-10)--(axis cs:0.490676743,-9.98795456)--(axis cs:0.980171403,-9.95184727)--(axis cs:1.46730474,-9.8917651)--(axis cs:1.95090322,-9.8078528)--(axis cs:2.4298018,-9.70031253)--(axis cs:2.90284677,-9.56940336)--(axis cs:3.36889853,-9.41544065)--(axis cs:3.82683432,-9.23879533)--(axis cs:4.27555093,-9.03989293)--(axis cs:4.71396737,-8.81921264)--(axis cs:5.14102744,-8.5772861)--(axis cs:5.55570233,-8.31469612)--(axis cs:5.95699304,-8.03207531)--(axis cs:6.34393284,-7.73010453)--(axis cs:6.71558955,-7.40951125)--(axis cs:7.07106781,-7.07106781)--(axis cs:7.40951125,-6.71558955)--(axis cs:7.73010453,-6.34393284)--(axis cs:8.03207531,-5.95699304)--(axis cs:8.31469612,-5.55570233)--(axis cs:8.5772861,-5.14102744)--(axis cs:8.81921264,-4.71396737)--(axis cs:9.03989293,-4.27555093)--(axis cs:9.23879533,-3.82683432)--(axis cs:9.41544065,-3.36889853)--(axis cs:9.56940336,-2.90284677)--(axis cs:9.70031253,-2.4298018)--(axis cs:9.8078528,-1.95090322)--(axis cs:9.8917651,-1.46730474)--(axis cs:9.95184727,-0.980171403)--(axis cs:9.98795456,-0.490676743)--(axis cs:10,0)--(axis cs:9.98795456,0.490676743)--(axis cs:9.95184727,0.980171403)--(axis cs:9.8917651,1.46730474)--(axis cs:9.8078528,1.95090322)--(axis cs:9.70031253,2.4298018)--(axis cs:9.56940336,2.90284677)--(axis cs:9.41544065,3.36889853)--(axis cs:9.23879533,3.82683432)--(axis cs:9.03989293,4.27555093)--(axis cs:8.81921264,4.71396737)--(axis cs:8.5772861,5.14102744)--(axis cs:8.31469612,5.55570233)--(axis cs:8.03207531,5.95699304)--(axis cs:7.73010453,6.34393284)--(axis cs:7.40951125,6.71558955)--(axis cs:7.07106781,7.07106781)--(axis cs:6.71558955,7.40951125)--(axis cs:6.34393284,7.73010453)--(axis cs:5.95699304,8.03207531)--(axis cs:5.55570233,8.31469612)--(axis cs:5.14102744,8.5772861)--(axis cs:4.71396737,8.81921264)--(axis cs:4.27555093,9.03989293)--(axis cs:3.82683432,9.23879533)--(axis cs:3.36889853,9.41544065)--(axis cs:2.90284677,9.56940336)--(axis cs:2.4298018,9.70031253)--(axis cs:1.95090322,9.8078528)--(axis cs:1.46730474,9.8917651)--(axis cs:0.980171403,9.95184727)--(axis cs:0.490676743,9.98795456)--(axis cs:6.123234e-16,10)--(axis cs:-0.490676743,9.98795456)--(axis cs:-0.980171403,9.95184727)--(axis cs:-1.46730474,9.8917651)--(axis cs:-1.95090322,9.8078528)--(axis cs:-2.4298018,9.70031253)--(axis cs:-2.90284677,9.56940336)--(axis cs:-3.36889853,9.41544065)--(axis cs:-3.82683432,9.23879533)--(axis cs:-4.27555093,9.03989293)--(axis cs:-4.71396737,8.81921264)--(axis cs:-5.14102744,8.5772861)--(axis cs:-5.55570233,8.31469612)--(axis cs:-5.95699304,8.03207531)--(axis cs:-6.34393284,7.73010453)--(axis cs:-6.71558955,7.40951125)--(axis cs:-7.07106781,7.07106781)--(axis cs:-7.40951125,6.71558955)--(axis cs:-7.73010453,6.34393284)--(axis cs:-8.03207531,5.95699304)--(axis cs:-8.31469612,5.55570233)--(axis cs:-8.5772861,5.14102744)--(axis cs:-8.81921264,4.71396737)--(axis cs:-9.03989293,4.27555093)--(axis cs:-9.23879533,3.82683432)--(axis cs:-9.41544065,3.36889853)--(axis cs:-9.56940336,2.90284677)--(axis cs:-9.70031253,2.4298018)--(axis cs:-9.8078528,1.95090322)--(axis cs:-9.8917651,1.46730474)--(axis cs:-9.95184727,0.980171403)--(axis cs:-9.98795456,0.490676743)--(axis cs:-10,1.2246468e-15)--cycle;
\path[fill=white,draw=black,line width=.5pt] (axis cs:9.9993155266,-.117)--(axis cs:19.8,-.117)--(axis cs:19.8,.117)--(axis cs:9.9993155266,.117)--cycle;
\end{axis}
\end{tikzpicture}
\caption{L2 current near the slot with the coil centred at $(15.2,0)$~mm. Samples are cell-centre values in the upper layer, at $z\simeq-0.333$~mm. Colour gives the magnitude of the real current vector, and arrows show its in-plane components with one common scale.}
\label{fig:ect.current}
\end{figure}
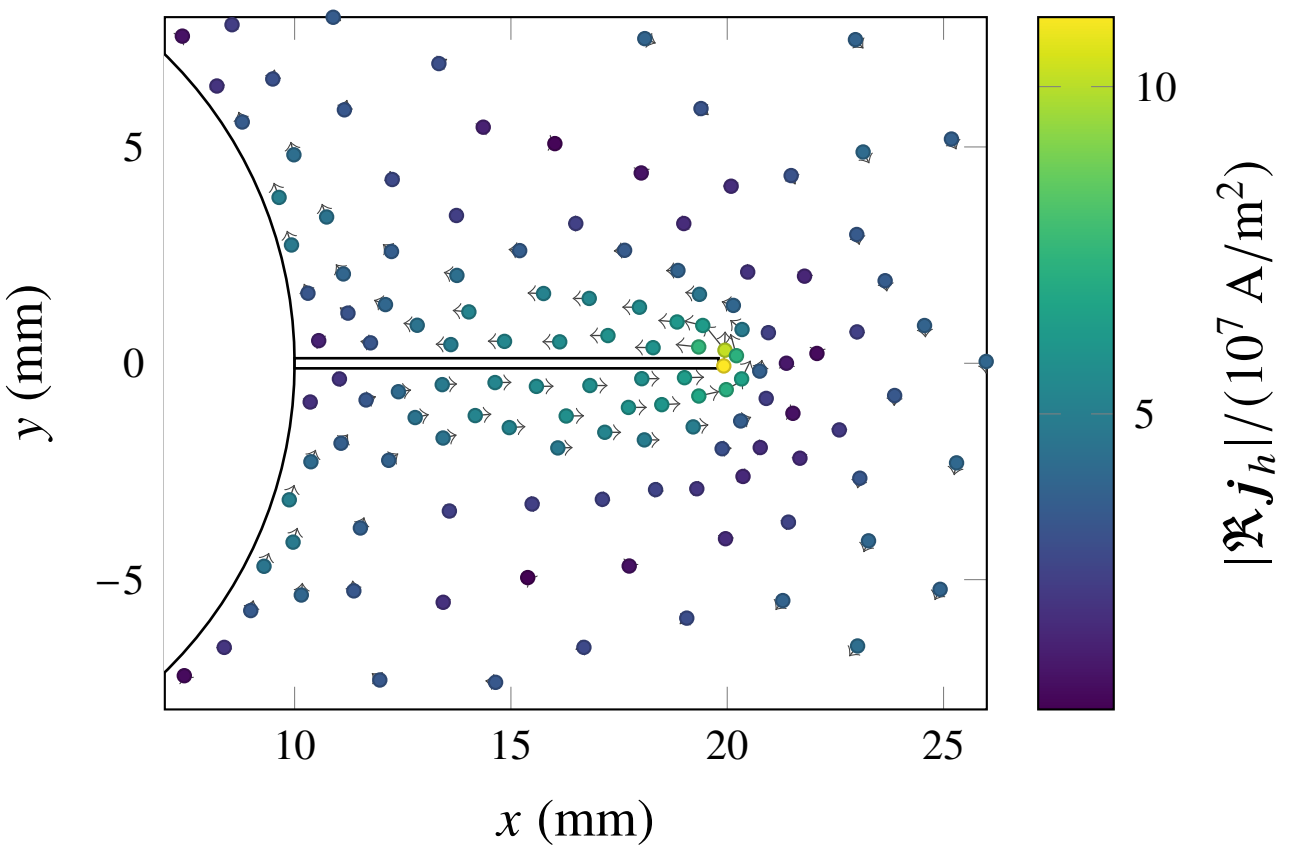